\documentclass[number,citesort,dvips]{arxbj}

\aid{0}
\volume{16}
\issue{2}
\pubyear{2010}
\firstpage{389}
\lastpage{417}
\doi{10.3150/09-BEJ223}

\makeatletter

\newtheorem{Theo}{Theorem}
\newtheorem{Prop}{Proposition}
\newtheorem{Lemma}{Lemma}

\newremark{example}{Example}
\newremark{remark}{Remark}

\makeatother

\begin{document}
\begin{frontmatter}

\title{Approximating a geometric fractional Brownian motion and related
processes via discrete Wick calculus}
\runtitle{Approximating a GFBM and related processes}

\begin{aug}
\author{\fnms{Christian} \snm{Bender}\ead[label=e1,mark]{bender@math.uni-sb.de}\thanksref{e1}}
\and
\author{\fnms{Peter} \snm{Parczewski}\ead[label=e2,mark]{parcz@math.uni-sb.de}\thanksref{e2}\corref{}}
\address{Department of Mathematics, Saarland University, P.O. Box 15
11 50, D-66041 Saarbr\"{u}cken, Germany. \printead{e1,e2}}

\runauthor{C. Bender and P. Parczewski}
\end{aug}

\received{\smonth{1} \syear{2009}}
\revised{\smonth{6} \syear{2009}}

%
\begin{abstract}
We approximate the solution of some linear systems of SDEs driven by a
fractional Brownian motion $B^{H}$ with Hurst
parameter $H \in(\frac{1}{2},1)$ in the Wick--It{\^{o}} sense,
including a geometric fractional Brownian motion. To this end, we apply
a Donsker-type approximation of
the fractional Brownian motion by disturbed binary random walks due
to Sottinen. Moreover, we replace the rather complicated Wick products
by their discrete counterpart, acting on the binary variables, in the
corresponding systems of Wick difference equations. As the solutions of
the SDEs admit series representations in terms of Wick powers, a key to
the proof of our Euler scheme is an approximation of the Hermite
recursion formula for the Wick powers of $B^{H}$.
\end{abstract}

%

\pdfkeywords{discrete Wick calculus, fractional Brownian motion, weak convergence, Wick--Ito integral}

\end{frontmatter}

\section{Introduction}

A fractional Brownian motion $B^{H}$ with Hurst parameter $H \in(0,1)$
is a continuous zero-mean Gaussian process in $\mathbb{R}$ with
stationary increments and covariance function
\[
\mathbf{E}[B^{H}_{t}B^{H}_{s}] = \tfrac{1}{2} (|t|^{2H} + |s|^{2H} -
|t-s|^{2H} ).
\]
The process $B^{1/2}$ is a standard Brownian motion, but a
fractional Brownian motion is not a semimartingale for $H \neq
\frac{1}{2}$. In this paper, we restrict ourselves to the case
$H>1/2$, in which the corresponding fractional Gaussian noise
$ (B^{H}_{n+1} - B^{H}_{n} )_{n \in\mathbb{N}}$ exhibits
long-range dependence.

In recent years, a lively interest in integration theory with respect
to fractional Brownian motion has emerged (see, e.g., the monographs by
Mishura or Biagini \textit{et al.} \cite{Mishura,Biagini}). One of the
extensions of the It{\^{o}} integral beyond semimartingales is the
fractional Wick--It{\^{o}} integral. It is based on the Wick product
$\diamond$, which has its origin as a renormalization operator in
quantum physics.
In probability theory, the Wick product with ordinary differentiation
rule imitates the situation of ordinary multiplication with It{\^{o}}
differentiation rule (cf.~Holden \textit{et al.} \cite{HoldenBuch}).
Actually, this makes it a natural tool to apply for extending the It{\^
{o}} integral.

We first consider the fractional Dol{\'{e}}ans--Dade SDE $\mathrm{d}S_{t} =
S_{t}\mathrm{d}^{\diamond}B^{H}_{t}$, $S_{0} = 1$, in terms of the fractional
Wick--It{\^{o}} integral. The well-known solution, $\exp(B^{H}_{t}
- \frac{1}{2}t^{2H} )$, is the geometric fractional Brownian
motion, also known as the Wick exponential of fractional Brownian
motion. Note that the Wick exponential has expectation equal to
one and can therefore be interpreted as a multiplicative noise.
Moreover, the ordinary exponential can be obtained from the Wick
exponential by a deterministic scaling. Neither process is a
semimartingale for $H \neq\frac{1}{2}$. The name ``Wick exponential''
is justified by the fact that it exhibits a power series expansion
with Wick powers $ (B^{H}_{t} )^{\diamond k}$ instead of
ordinary powers.

More generally, we consider a linear system of SDEs,
\begin{eqnarray}\label{lineare_system}
\mathrm{d}X_{t} &=& (A_{1}X_{t} + A_{2}Y_{t} )\,\mathrm{d}^{\diamond}B^{H}_{t} ,\qquad X_{0} = x_{0},\nonumber\\[-8pt]\\[-8pt]
\mathrm{d}Y_{t} &=& (B_{1}X_{t} + B_{2}Y_{t} )\,\mathrm{d}^{\diamond}B^{H}_{t} ,\qquad Y_{0} = y_{0}.\nonumber
\end{eqnarray}
One can obtain Wick power series expansions for the solution of this
system, too. Our goal is to approximate these Wick analytic
functionals of a fractional Brownian motion. To this end, we require
an approximation of a fractional Brownian motion and an
approximation of the Wick product.

There are several ways to approximate a fractional Brownian motion. One
of the first approximations was given by Taqqu \cite{Taqqu} in terms
of stationary Gaussian sequences. We refer to Mishura \cite{Mishura}, Section
1.15.3, for further approaches to weak convergence to a
fractional Brownian motion. Sottinen constructed a simple approximation
of a fractional Brownian motion on an interval for $H > \frac{1}{2}$ by
sums of square-integrable random variables in \cite{So}. He used the
Wiener integral representation of a fractional Brownian motion on an interval,
$B^{H}_{t} = \int_{0}^{t}z_{H}(t,s)\,\mathrm{d}B_{s}$,
for a suitable deterministic kernel $z_H(t,s)$, due to Molchan and
Golosov, and Norros \textit{et al.} \cite{Molchan,MolchanGolosov,Norros}.
For this purpose, he combined a pointwise approximation of the
kernel $z_{H}(t,s)$ with Donsker's theorem. This approach was
extended by Nieminen \cite{Nieminen} to weak convergence of
perturbed martingale differences to fractional Brownian motion. We
shall utilize Sottinen's approximation with binary random variables
throughout this paper.

The main problem of applying the Wick product on random variables with
continuous distributions is that it is not a pointwise operation.
Thus, an explicit computation of the Wick--It{\^{o}} integral is only
possible in rare special cases. But this is precisely the advantage of
the binary random walks.
In such a purely discrete setup, we apply the discrete counterpart of the
Wick product as introduced in Holden \textit{et al.} \cite{HoldenSDE}.
Starting from the binary random walk, one can build up a discrete
Wiener space, and the discrete Wick product depends on this
discretization. This Wiener chaos gives the analogy to the
continuous Wick products. For a survey on discrete Wiener chaos, we
refer to Gzyl \cite{Gzyl}. However, we will introduce the discrete
Wick product in a self-contained way in Section \ref{approxsection}.

We can now formulate a weak Euler scheme of the linear system of SDEs
(\ref{lineare_system}) in the Wick--It{\^{o}} sense,
\begin{eqnarray}\label{linearesystemWick}
X^{n}_{l} &=& X^{n}_{l-1} + (A_{1}X^{n}_{l-1} + A_{2}Y^{n}_{l-1} )
\diamond_{n} \bigl(B^{H,n}_{l/n} - B^{H,n}_{(l-1)/n} \bigr)
,\nonumber\\
X^{n}_{0} &=& x_{0},\qquad l = 1, \ldots, n,\nonumber\\[-8pt]\\[-8pt]
Y^{n}_{l} &=& Y^{n}_{l-1} + (B_{1}X^{n}_{l-1} + B_{2}Y^{n}_{l-1} )
\diamond_{n} \bigl(B^{H,n}_{l/n} - B^{H,n}_{(l-1)/n} \bigr)
,\nonumber\\
Y^{n}_{0} &=& y_{0},\qquad l = 1, \ldots, n,\nonumber
\end{eqnarray}
where $\diamond_{n}$ is the discrete Wick product and
$ (B^{H,n}_{l/n} - B^{H,n}_{(l-1)/n} )$ are
the increments of the disturbed binary random walk. As a main result,
we show that the piecewise constant interpolation of the solution of
(\ref{linearesystemWick}) converges weakly in the Skorokhod space
to the solution of (\ref{lineare_system}). This is the first
rigorous convergence result connecting discrete and continuous Wick
calculus of which we are aware. As a special case, (\ref{linearesystemWick}) contains the Wick difference equation
%
\begin{equation}\label{simpleWickdiffeq}
X^{n}_{l} = X^{n}_{l-1} + X^{n}_{l-1}\diamond_{n} \bigl(B^{H,n}_{l/n} - B^{H,n}_{(l-1)/n} \bigr) ,\qquad X^{n}_{0} = 1,\qquad l = 1, \ldots, n.
\end{equation}
As a consequence, the piecewise constant interpolation of (\ref{simpleWickdiffeq}) converges weakly to a geometric fractional Brownian
motion, the solution of the fractional Dol{\'{e}}ans--Dade SDE. This
was conjectured by Bender and Elliott \cite{BenderElliott} in their
study of the Wick fractional Black--Scholes market.

In \cite{So}, Sottinen considered the corresponding difference
equation in the pathwise sense, that is, with ordinary
multiplication instead of the discrete Wick product:
%
\begin{equation}\label{simpleWickdiffeqsottinen}
\hat{X}^{n}_{l} = \hat{X}^{n}_{l-1} + \hat{X}^{n}_{l-1} \bigl(B^{H,n}_{l/n} - B^{H,n}_{(l-1)/n} \bigr) ,\qquad \hat{X}^{n}_{0} = 1,\qquad l = 1, \ldots
, n.
\end{equation}
The solution is explicitly given by the multiplicative expression
%
\begin{equation}\label{Sottinendiffsolasprod}
\hat{X}_{l}^{n} = \prod_{j = 1}^{l} \bigl(1 + \bigl(B^{H,n}_{j/n} -
B^{H,n}_{(j-1)/n} \bigr) \bigr).
\end{equation}
By the logarithmic transform of ordinary products into sums and a
Taylor expansion, one obtains an additive expression for
$\ln(\hat{X}_{l}^{n})$ which converges weakly to a fractional
Brownian motion. In this way, Sottinen proved the convergence of
$\hat{X}$ to the ordinary exponential of a fractional Brownian
motion \cite{So}, Theorem 3. This approach fails for the solution of
(\ref{simpleWickdiffeq}) since, in a product representation,
analogous to (\ref{Sottinendiffsolasprod}), the discrete Wick
product $\diamond_{n}$ appears instead of ordinary
multiplication. There is, however, no straightforward way to
transform discrete Wick products into sums by application of a
continuous functional.

However, the solution of (\ref{linearesystemWick}) exhibits an
expression which is closely related to a discrete Wick power series
representation. Therefore, the convergence can be initiated
explicitly for the Wick powers of a fractional Brownian motion,
which fulfill the Hermite recursion formula. We obtain a discrete
analog to this recursion formula for discrete Wick powers of
disturbed binary random walks. Actually, the weak convergence of
these discrete Wick powers is the key to the proof for our Euler
scheme.

The paper is organized as follows. In Section \ref{Wickexpsection},
we give some preliminaries on the Wick--It{\^{o}} integral with respect
to a fractional Brownian motion and introduce the Wick exponential and
other Wick analytic functionals.
We then define the approximating sequences and state the main results
in Section \ref{approxsection}. Section \ref{WickDiffsec}
is devoted to some $L^2$- estimates of the approximating sequences. We
prove convergence in finite-dimensional distributions in Section \ref{secconvfd} and tightness in Section \ref{sectightness}.

\section{Wick exponential and Wick analytic functionals}\label{Wickexpsection}

In this section, we introduce the Wick product and the Wick--It{\^{o}}
integral, and describe the Hermite recursion formula for Wick powers of
a zero-mean Gaussian random variable. We then obtain the Wick power
series expansions for the solutions of SDEs (\ref{lineare_system}).

We consider a \textit{geometric fractional Brownian motion} or the
so-called \textit{Wick exponential} of a fractional Brownian motion
$\exp(B^{H}_{t} - \frac{1}{2}t^{2H} )$.
For $H = \frac{1}{2}$, this is exactly a geometric Brownian motion,
also known as the stochastic exponential of a standard Brownian
motion. For all $H \in(0,1)$ and $t \geq0$, it holds that $t^{2H}
= \mathbf{E} [ (B^{H}_{t} )^2 ]$ and thus the Wick
exponential generalizes the stochastic exponential. It is well known
that $\exp(B_{t} - \frac{1}{2}t )$ solves the
Dol{\'{e}}ans--Dade equation
\[
\mathrm{d}S_{t} = S_{t}\,\mathrm{d}B_{t} ,\qquad S_{0} = 1,
\]
where the integral is an ordinary It{\^{o}} integral. Actually, the
Wick exponential of fractional Brownian motion solves the corresponding
fractional Dol{\'{e}}ans--Dade equation
\[
\mathrm{d}S_{t} = S_{t}\,\mathrm{d}^{\diamond}B^{H}_{t} ,\qquad S_{0} = 1,
\]
in terms of a fractional Wick--It{\^{o}} integral (cf.~Mishura
\cite{Mishura}, Theorem 3.3.2). We want to approximate solutions of similar SDEs.

Let $\Phi$ and $\Psi$ be two zero-mean Gaussian random variables. The
\textit{Wick exponential} is then defined as
\[
\exp^{\diamond}(\Phi) := \exp\bigl(\Phi- \tfrac{1}{2}\mathbf{E}[|\Phi|^2] \bigr).
\]
For a standard Brownian motion $(B_{t})_{t \geq0}$ and $s<t<u$, it
holds that
\[
\exp^{\diamond}(B_{u} - B_{t}) \exp^{\diamond}(B_{t} - B_{s}) = \exp
^{\diamond}(B_{u}- B_{s}).
\]
Forcing this renormalization property to hold for all, possibly
correlated, $\Phi$ and $\Psi$, leads to the definition of the \textit
{Wick product} $\diamond$ of two Wick exponentials:
\[
\exp^{\diamond}(\Phi) \diamond\exp^{\diamond}(\Psi) := \exp^{\diamond
}(\Phi+ \Psi).
\]
The Wick product can be extended to larger classes of random variables
by density arguments (cf.~\cite{Duncan,Nualart,Bender}). For a
general introduction to the Wick product, we refer to the monographs by
Kuo and Holden \textit{et al.} \cite{Kuo,HoldenBuch} and Hu and Yan
\cite{HuYan}. Note that the Wick product is not a pointwise
operation. If we suppose that $\Phi\sim\mathcal{N}(0,\sigma)$, then
we have, by definition, $\Phi^{\diamond0} = 1$, $\Phi^{\diamond1} =
\Phi$ and the recursion
\[
\Phi^{\diamond n+1} = \Phi^{\diamond n} \diamond\Phi.
\]
Observe that it holds that
\begin{eqnarray*}
\frac{\mathrm{d}}{\mathrm{d}x} \exp^{\diamond} (x \Phi) \bigg|_{x = 0} &=& \frac{\mathrm{d}}{\mathrm{d}x} \exp
\biggl(x \Phi- \frac{1}{2}\mathbf{E} [|x\Phi|^2 ] \biggr) \bigg|_{x = 0} \\
&=& (\Phi- x\sigma^2 )\exp\biggl(x \Phi- \frac{1}{2}\mathbf{E} [|x\Phi|^2
] \biggr) \bigg|_{x = 0} =\Phi.
\end{eqnarray*}
Suppose we have
\[
\Phi^{\diamond k} = \frac{\mathrm{d}^{k}}{\mathrm{d}w^{k}}\exp^{\diamond} (w\Phi) \bigg|_{w=0}
\]
for all positive integers $k \leq n$. Then, with $z = w+x$, $\frac
{\mathrm{d}z}{\mathrm{d}w}=\frac{\mathrm{d}z}{\mathrm{d}x} = 1$, we get
\begin{eqnarray*}
\Phi^{\diamond(n+1)} &=& \frac{\mathrm{d}^{n}}{\mathrm{d}w^{n}}\exp^{\diamond} (w\Phi)
\bigg|_{w=0} \diamond \frac{\mathrm{d}}{\mathrm{d}x} \exp^{\diamond} (x \Phi) \bigg|_{x = 0}\\
&=& \frac{\mathrm{d}^{n}}{\mathrm{d}w^{n}}\frac{\mathrm{d}}{\mathrm{d}x} \exp^{\diamond} \bigl((w+x)\Phi\bigr)
\bigg|_{w=0, x=0} = \frac{\mathrm{d}^{n+1}}{\mathrm{d}z^{n+1}}\exp^{\diamond} (z\Phi) \bigg|_{z=0}.
\end{eqnarray*}
We now obtain, by differentiation and the Leibniz rule, the following
Wick recursion formula:
\begin{eqnarray}\label{recursiononwickpower}
\Phi^{\diamond n+1} &=& \frac{\mathrm{d}^n}{\mathrm{d}w^n} \biggl( (\Phi-w\sigma^2 )\exp
\biggl(w\Phi- \frac{1}{2}\mathbf{E} [|w\Phi|^2 ] \biggr) \biggr) \bigg|_{w = 0} \nonumber\\
&=& (\Phi-w\sigma^2 )\frac{\mathrm{d}^n}{\mathrm{d}w^n}\exp\biggl(w\Phi- \frac{1}{2}\mathbf
{E} [|x\Phi|^2 ] \biggr) \bigg|_{w = 0}\nonumber\\[-8pt]\\[-8pt]
&&{} + n(-\sigma^2)\frac{\mathrm{d}^{n-1}}{\mathrm{d}w^{n-1}} \exp\biggl(w\Phi- \frac
{1}{2}\mathbf{E} [|w\Phi|^2 ] \biggr) \bigg|_{w = 0}\nonumber\\
&=& \Phi\Phi^{\diamond n} - n\sigma^2\Phi^{\diamond n-1}.\nonumber
\end{eqnarray}

Define the \textit{Hermite polynomial of degree $n \in\mathbb{N}$ with
parameter $\sigma^2$} as
\[
h^{n}_{\sigma^{2}}(x) := (-\sigma^2)^{n} \exp\biggl(\frac{x^2}{2\sigma^2} \biggr)
\frac{\mathrm{d}^n}{\mathrm{d}x^n} \exp\biggl(\frac{- x^2}{2\sigma^2} \biggr).
\]
The series expansion
%
\begin{equation}\label{wickexpsoalsreihe}
\exp\biggl(x - \frac{1}{2}\sigma^{2}\biggr) = \sum_{n =
0}^{\infty}\frac{1}{n!} h^{n}_{\sigma^{2}}(x)
\end{equation}
then holds true. The first Hermite polynomials are $h^{0}_{\sigma^{2}}(x)
= 1$, $h^{1}_{\sigma^2}(x) = x$. By the Leibniz rule, we obtain the
Hermite recursion formula
%
\begin{equation}\label{hermiteoriginrecursion}
h^{n+1}_{\sigma^{2}}(x) = x h^{n}_{\sigma^{2}}(x) - n \sigma^2
h^{n-1}_{\sigma^{2}}(x).
\end{equation}

By the equivalent first terms and recursions (\ref{recursiononwickpower}) and (\ref{hermiteoriginrecursion}), we can conclude that for
any Gaussian random variable $\Phi\sim\mathcal{N}(0,\sigma)$ and all
$n \in\mathbb{N}$,
we have
%
\begin{equation}\label{propwickpowersarehermitepolynomials}
\Phi^{\diamond n} = h^{n}_{\sigma^2}(\Phi).
\end{equation}
By (\ref{wickexpsoalsreihe}), we additionally have
%
\begin{equation}\label{wickexpalswickpowerseries}
\exp^{\diamond}(\Phi) = \sum_{n = 0}^{\infty}\frac{1}{n!}\Phi^{\diamond n}.
\end{equation}
The fractional Wick--It{\^{o}} integral, introduced by Duncan \textit
{et al.}
\cite{Duncan}, is an extension of the It{\^{o}} integral beyond the
semimartingale framework. There are several approaches to the
fractional Wick--It{\^{o}} integral. Essentially, these approaches are
via white noise theory, as in Elliott and von der Hoek \cite{ElliottHoek}, and Hu and \O ksendal \cite{HuOksendal}, by Malliavin
calculus in Al\`{o}s \textit{et al.} \cite{Alosetal}, or by an
S-transform approach
in Bender \cite{Bender}. In contrast to the forward integral, the
fractional Wick--It{\^{o}} integral has zero mean in general. This is
the crucial property for an additive noise. The Wick--It{\^{o}}
integral is based on the Wick product. For a sufficiently good
process $(X_{s})_{s \in[0,t]}$, the fractional Wick--It{\^{o}} integral
with respect to fractional Brownian motion $(B_{s}^{H})_{[0,t]}$ can
be easily defined by Wick--Riemann sums (cf.~Duncan \textit{et al.}
\cite{Duncan} or Mishura \cite{Mishura}, Theorem 2.3.10). If we
suppose that
$\pi_{n} = \{0 = t_{0} < t_{1} < \cdots< t_{n} = t\}$ with
$\max_{t_{i} \in\pi_{n}}|t_{i} - t_{i-1}| \rightarrow0$ for
$n \rightarrow\infty$, then
\[
\int_{0}^{t}X_{s}\,\mathrm{d}^{\diamond}B_{s}^{H} :=\lim_{n \rightarrow\infty}
\sum_{t_{i} \in\pi_{n}} X_{t_{i-1}} \diamond(B_{t_{i}}^{H} -
B_{t_{i-1}}^{H} ),
\]
if the Wick products and the $L^2(\Omega)$-limit exist. For more
information on Wick--It{\^{o}} integrals with respect to fractional
Brownian motion, we refer to Mishura \cite{Mishura}, Chapter 2.

By the fractional It{\^{o}} formula (cf.~\cite{Bender}, Theorem 5.3 or
\cite{Biagini}, Theorem 3.7.2), we have
%
\begin{equation}\label{Wickpowerchainrule}
\mathrm{d}(B^{H}_{t})^{\diamond k} =k(B^{H}_{t})^{\diamond k-1} \,\mathrm{d}^{\diamond}B^{H}_{t} ,\qquad
(B^{H}_{0})^{\diamond k} = \mathbf{1}_{ \{k = 0 \}}.
\end{equation}
For the Wick exponential
%
\begin{equation}\label{WickexpfBMalsreihe}
\exp^{\diamond} (B^{H}_{t} ) = \sum_{k=0}^{\infty}\frac{1}{k!}
(B^{H}_{t} )^{\diamond k},
\end{equation}
we obtain, by summing up the identity (\ref{Wickpowerchainrule}),
the fractional Dol{\'{e}}ans--Dade equation,
%
\begin{equation}\label{SDEforWickexponential}
\mathrm{d}S_{t} = S_{t}\,\mathrm{d}^{\diamond}B^{H}_{t} ,\qquad S_{0} = 1.
\end{equation}
For any analytic function $F(x) =\sum_{k = 0}^{\infty} \frac
{a_{k}}{k!}x^{k}$, we define the \textit{Wick version} as
\[
F^{\diamond}(\Phi) =\sum_{k = 0}^{\infty} \frac{a_{k}}{k!}\Phi
^{\diamond k}.
\]
From the recursive system of SDEs (\ref{Wickpowerchainrule}), we
obtain SDEs for other \textit{Wick analytic functionals} of a
fractional Brownian motion
\[
F^{\diamond}(B^{H}_{t}) =\sum_{k = 0}^{\infty} \frac{a_{k}}{k!}
(B^{H}_{t} )^{\diamond k} .
\]

Recall the linear system of SDEs (\ref{lineare_system}),
\begin{eqnarray*}
\mathrm{d}X_{t} &=& (A_{1}X_{t} + A_{2}Y_{t} )\,\mathrm{d}^{\diamond}B^{H}_{t} ,\qquad X_{0} = x_{0},\\
\mathrm{d}Y_{t} &=& (B_{1}X_{t} + B_{2}Y_{t} )\,\mathrm{d}^{\diamond}B^{H}_{t} ,\qquad Y_{0} = y_{0}.
\end{eqnarray*}

The coefficients of the solution,
%
\begin{equation}\label{linearSDEsystemsolution}
X_{t} = \sum_{k = 0}\frac{a_{k}}{k!} (B^{H}_{t} )^{\diamond k},
\qquad
Y_{t} = \sum_{k = 0}\frac{b_{k}}{k!} (B^{H}_{t} )^{\diamond k},
\end{equation}
can be obtained recursively via (\ref{Wickpowerchainrule}) to be
\[
a_{0} = x_{0} ,\qquad
b_{0} = y_{0} ,\qquad
a_{k} = A_{1}a_{k-1} + A_{2}b_{k-1} ,\qquad
b_{k} = B_{1}a_{k-1} + B_{2}b_{k-1}.
\]

Note that it holds that $|a_{k}| , |b_{k}| \leq C^{k}$ for a $C \in
\mathbb{R}_{+}$. This is according to the recursive derivation of
the coefficients and it ensures that the Wick analytic functionals
$X_{t}$ and $Y_{t}$ are square-integrable (cf.~the proof of
Proposition \ref{fdmaintheorem}).

\section{The approximation results}\label{approxsection}

Here, we present the approximating sequences and discuss the main
results. More precisely, we introduce the Donsker-type approximation of
a fractional Brownian motion and the discrete Wick product, and obtain
Wick difference equations, which correspond to the SDEs. We shall work
with a fractional Brownian motion on the interval $[0,1]$, but all
results extend to any compact interval $[0,T]$.

We first consider the following kernel representation of a fractional
Brownian motion on the interval $[0,1]$, based on works by Molchan and
Golosov \cite{Molchan,MolchanGolosov},
%
\begin{equation}\label{fBMalsWiener}
B^{H}_{t} = \int^{t}_{0}z_{H}(t,s)\,\mathrm{d}B_{s}.
\end{equation}
For $H > \frac{1}{2}$, the deterministic kernel takes the form
%
\begin{equation}\label{fBM_kernel}
z_{H}(t,s) = \mathbf{1}_{ \{t \geq s \}} c_{H}\biggl(H - \frac{1}{2}\biggr)s^{
1/2 - H} \int^{t}_{s}u^{H - 1/2}(u-s)^{H - 3/2}\,\mathrm{d}u
\end{equation}
with the constant
\[
c_{H} = \sqrt{\frac{2H\Gamma(3/2 - H)}{\Gamma(H + 1/2)\Gamma(2 - 2H)}},
\]
where $\Gamma$ is the Gamma function (Norros \textit{et al.} \cite{Norros} or
Nualart \cite{Nualart}, Section 5.1.3). In order to simplify the
notation, we think of $H \in(\frac{1}{2},1)$ as fixed from now on
and omit the subscript $H$ in the notation of the kernel. For an
introduction to some elementary properties of fractional Brownian
motion, we refer to Nualart \cite{Nualart}, Chapter 5, Mishura
\cite{Mishura} or Biagini \textit{et al.} \cite{Biagini}.

We apply Sottinen's approximation of a fractional Brownian motion by
\textit{disturbed binary random walks}. Suppose $(\Omega, \mathcal{F},
P)$ is a probability space and, for all $n \in\mathbb{N}$ and $i =
1,\ldots, n$, we have independent and identically distributed symmetric
Bernoulli random variables $\xi^{n}_{i}\dvtx \Omega\rightarrow\{-1,1\}$
with $P(\xi^{n}_{i} = 1) = P(\xi^{n}_{i} = -1)$. By Donsker's theorem,
the sequence of random walks
$B^{(n)}_{t} = \frac{1}{\sqrt{n}}\sum^{ \lfloor nt \rfloor}_{i = 1}\xi
^{(n)}_{i}$
converges weakly to a standard Brownian motion $B = (B_{t})_{t \in
[0,1]}$ \cite{Billialt}, Theorem~16.1. The idea of Sottinen \cite{So} is to combine these random walks with a pointwise approximation of
the kernel in representation (\ref{fBMalsWiener}). Define the
pointwise approximation of $z(t,s)$ as
\[
z^{(n)}(t,s) := n\int^{s}_{s-1/n}z \biggl(\frac{ \lfloor nt \rfloor
}{n}, u \biggr)\,\mathrm{d}u.
\]
The sequence of binary random walks
\[
B^{H,n}_{t} := \int^{t}_{0}z^{(n)}(t,s)\,\mathrm{d}B^{(n)}_{s} = \sum^{ \lfloor nt
\rfloor}_{i = 1} n \int^{i/n}_{(i-1)/n}z \biggl(\frac{ \lfloor
nt \rfloor}{n}, s \biggr)\,\mathrm{d}s \frac{1}{\sqrt{n}} \xi^{(n)}_{i}
\]
then converges weakly to a fractional Brownian motion $(B_{t}^{H})_{t
\in[0,1]}$ in the Skorokhod space $D([0,1],\mathbb{R})$ \cite{So}, Theorem 1.

A major advantage of the binary random walks is that we can avoid the
difficult Wick product for random variables with continuous
distributions. We approximate this operator on the binary random walks
by \textit{discrete Wick products}.

For any $n \in\mathbb{N}$, let $ (\xi_{1}^{n}, \xi_{2}^{n}, \ldots, \xi
_{n}^{n} )$ be
the $n$-tuple of independent and identically distributed symmetric
Bernoulli random variables for the binary random walk $B^{H,n}_{t}$.
The \textit{discrete Wick product} is defined as
\[
\prod_{i \in A} \xi^{n}_{i} \diamond_{n} \prod_{i \in B} \xi^{n}_{i} :=
\cases{
\displaystyle\prod_{i \in A \cup B} \xi^{n}_{i}, &\quad if $ A \cap B = \emptyset$,\cr
0, &\quad otherwise,
}
\]
where $A,B \subseteq \{ 1,\ldots, n \}$. We denote by
\[
\mathcal{F}_{n} := \sigma(\xi_{1}^{n}, \xi_{2}^{n}, \ldots, \xi
_{n}^{n} )
\]
the $\sigma$-field generated by the Bernoulli variables. Define
\[
\Xi_{A}^n := \prod_{i \in A} \xi^{n}_{i}.
\]
Clearly, the family of functions $ \{\Xi_{A}^n \dvtx A \subseteq\{1,\ldots
, n\} \}$
is an orthonormal set in $L^2(\Omega, \mathcal{F}_{n},P)$. Since its
cardinality is equal to the dimension of $L^2(\Omega, \mathcal
{F}_{n},P)$, it constitutes a basis. Thus, every $X \in L^2(\Omega,
\mathcal{F}_{n},P)$ has a unique expansion, called the \textit{Walsh
decomposition},
\[
X = \sum_{A \subseteq\{1,\ldots, n\}} x_{A}^n\Xi_{A}^n,
\]
where $x^{n}_{A} \in\mathbb{R}$. The Walsh decomposition can be
regarded as a discrete version of the chaos expansion.
By algebraic rules, one obtains, for $X = \sum_{A \subseteq\{1,\ldots,
n\}} x_{A}^n\Xi_{A}^n$ and $Y = \sum_{B \subseteq\{1,\ldots, n\}
}y_{B}^n\Xi_{B}^n$,
\[
X \diamond_{n} Y = \sum_{C \subseteq\{1,\ldots, n\}} \Biggl(\mathop{\sum_{
A \cup B = C }}_{ A \cap B = \emptyset}x^{n}_{A}y^{n}_{B} \Biggr)\Xi_{C}^n.
\]
Furthermore, the $L^2$-inner product can be computed in terms of the
Walsh decomposition as
%
\begin{equation}\label{L2norminWalshterms}
\mathbf{E} [XY ] = \sum_{A \subseteq\{1,\ldots, n\}}x_{A}^n y_{A}^{n}.
\end{equation}
There exists an analogous formula for the Wick product on the white
noise space via chaos expansions that justifies the analogy between the
discrete and ordinary Wick calculus (cf.~Kuo \cite{Kuo}). For more
information on the discrete Wick product, we refer to Holden \textit{et
al.} \cite{HoldenSDE}. More generally, the introduction of a discrete
Wiener chaos depends on the class of discrete random variables $ (\xi
_{1}^{n}, \xi_{2}^{n}, \ldots, \xi_{n}^{n} )$. We refer to Gzyl \cite{Gzyl} for a survey of other discrete Wiener chaos approaches.

The representation
\[
B^{H,n}_{t} =\sum^{ \lfloor nt \rfloor}_{i = 1} b_{t,i}^{n} \xi^{n}_{i}
\qquad\mbox{with } b_{t,i}^{n} := \sqrt{n} \int^{i/n}_{
(i-1)/n}z\biggl(\frac{ \lfloor nt \rfloor}{n}, s\biggr)\,\mathrm{d}s
\]
is the Walsh decomposition for the binary random walk approximating
$B^{H}$ in $L^2(\Omega, \mathcal{F}_{n},P)$. Note that $b_{t,i}^{n} =
b_{\lfloor nt \rfloor/n,i}^{n}$. Thus, we can consider
$B_{t}^{H,n} = B_{\lfloor nt \rfloor/n}^{H,n}$
as a process in discrete time. We can now state our first convergence result.

\begin{Theo}\label{dattheo1}
Suppose
that:
\begin{enumerate}
\item$\lim_{n \rightarrow\infty}a_{n,k} = a_{k}$ exists for all $k
\in\mathbb{N}$;
\item there exists a $C \in\mathbb{R}_{+}$ such that $|a_{n,k}| \leq
C^{k}$ for all $n,k \in\mathbb{N}$.
\end{enumerate}
The sequence of processes
$\sum_{k=0}^{n}\frac{a_{n,k}}{k!} (B^{H,n} )^{\diamond_{n} k}$
then converges weakly to the Wick power series
$\sum_{k=0}^{\infty}\frac{a_{k}}{k!} (B^{H} )^{\diamond k}$
in the Skorokhod space $D([0,1],\mathbb{R})$.
\end{Theo}

The proof is given in Sections \ref{secconvfd} and \ref{sectightness}.

Consider now the following recursive system of Wick difference equations:
%
\begin{equation}\label{recursivesystemwickdifferenceforWickpowers}
U_{l}^{k,n} = U_{l-1}^{k,n} + k U_{l-1}^{k-1,n}\diamond_{n} \bigl(B_{l/n}^{H,n} - B_{(l-1)/n}^{H,n} \bigr) ,\qquad U^{0,n}_{l} = 1
,\qquad U^{k,n}_{0} = 0,
\end{equation}
for all $l = 1, \ldots, n$ and $k \in\mathbb{N}$. This is the discrete
counterpart of the recursive system of SDEs in (\ref{Wickpowerchainrule}). We observe that $U^{0,n} = 1 = (B^{H,n} )^{\diamond_{n} 0}$ and
$U^{1,n} = (B^{H,n} )^{\diamond_{n} 1}$, but
\[
U^{2,n}_{2} =2B_{1/n}^{H,n}\diamond_{n} B_{2/n}^{H,n}
\neq B_{2/n}^{H,n} \diamond_{n}B_{2/n}^{H,n} =
(B_{2/n}^{H,n} )^{\diamond_{n} 2}.
\]
Thus, in contrast to the continuous case in (\ref{Wickpowerchainrule}), the discrete Wick powers
are not the solutions for (\ref{recursivesystemwickdifferenceforWickpowers}) if $k \geq2$.

However, we can prove a variant of Theorem \ref{dattheo1}, based on
the system of recursive Wick difference equations, whose proof will
also be given in Sections \ref{secconvfd} and \ref{sectightness}.

\begin{Theo}\label{extendedmaintheorem}
Under the assumptions of Theorem \ref{dattheo1}, define $\widetilde
{U}^{k,n}_{t} := U^{k,n}_{ \lfloor nt \rfloor}$ as the piecewise
constant interpolation of (\ref{recursivesystemwickdifferenceforWickpowers}).

The sequence of processes
$\sum_{k=0}^{n}\frac{a_{n,k}}{k!} \widetilde{U}^{k,n}$
then converges weakly to the Wick power series
$\sum_{k=0}^{\infty}\frac{a_{k}}{k!} (B^{H} )^{\diamond k}$
in the Skorokhod space $D([0,1],\mathbb{R})$.
\end{Theo}

\begin{example}[(Wick powers of a fractional Brownian motion)]
For $a_{n,k} = l! \mathbf{1}_{ \{k = l \}}$,
\begin{eqnarray*}
(B^{H,n} )^{\diamond_{n} l} &\stackrel{d}{\longrightarrow}& (B^{H}
)^{\diamond l},\\
\widetilde{U}^{l,n} &\stackrel{d}{\longrightarrow}& (B^{H} )^{\diamond l}.
\end{eqnarray*}
\end{example}

\begin{example}[(Geometric fractional Brownian motion)]
For $a_{n,k} = a_{k} = 1$, we have
\begin{eqnarray*}
\exp^{\diamond_{n}} (B^{H,n}_{t} ) &:=& \sum_{k = 0}^{ \lfloor nt \rfloor
}\frac{1}{k!} (B^{H,n}_{t} )^{\diamond_{n} k} \stackrel{d}{\rightarrow
} \exp^{\diamond} (B^{H} ),\\
\widetilde{S}^{n} &:=& \sum_{k=0}^{n}\frac{1}{k!}\widetilde{U}^{k,n}
\stackrel{d}{\longrightarrow} \exp^{\diamond} (B^{H} ).
\end{eqnarray*}
Observe that by summing up the recursive system of Wick difference
equations (\ref{recursivesystemwickdifferenceforWickpowers}), we obtain
%
\begin{equation}\label{thediffeq}
S_{l}^{n} = S_{l-1}^{n} + S_{l-1}^{n}\diamond_{n} \bigl(B_{
l/n}^{H,n} - B_{(l-1)/n}^{H,n} \bigr),\qquad S_{0}^{n} = 1,
\end{equation}
for $l = 1, \ldots, n$, where $S^{n}_{l} = \widetilde{S}^{n}_{l/n}$.
Hence, the piecewise constant interpolation of (\ref{thediffeq})
converges weakly to the solution of the fractional Dol{\'{e}}ans--Dade
equation (\ref{SDEforWickexponential}).
\end{example}

The reasoning of the previous example can be generalized as follows.

\begin{Theo}[(Linear SDE with drift)]\label{approxlinearSDEwithdrift}
Suppose $\mu, s_{0} \in\mathbb{R}$, $\sigma> 0$. Then $\widetilde
{S}^{n}_{t} := S^{n}_{ \lfloor nt \rfloor}$, where $S^{n}$ is the
solution of the Wick difference equation
%
\begin{equation}\label{discretelinearWickdifferencedrift}
S_{l}^{n} = \biggl(1+\frac{\mu}{n} \biggr) S_{l-1}^{n} + \sigma S_{l-1}^{n}\diamond
_{n} \bigl(B_{l/n}^{H,n} - B_{(l-1)/n}^{H,n} \bigr),\qquad S_{0}^{n} =
s_{0},\qquad l = 1,\ldots, n,
\end{equation}
converges weakly to the solution of the \textit{linear SDE with drift}
%
\begin{equation}\label{wicklinearSDEdrift}
\mathrm{d}S_{t} = \mu S_{t} \,\mathrm{d}t + \sigma S_{t} \,\mathrm{d}^{\diamond}B^{H}_{t},\qquad S_{0} = s_{0},
\end{equation}
in the Skorokhod space $D([0,1],\mathbb{R})$.
\end{Theo}

\begin{pf}
First, observe that for $\sigma_{n} \rightarrow\sigma> 0$ and
$\widetilde{a}_{n,k} = a_{n,k}\sigma_{n}^{k}$, we obtain, by Theorem
\ref{extendedmaintheorem}, that
\[
\widetilde{V}^{n} := \sum_{k=0}^{n}\frac{a_{n,k}}{k!} \sigma_{n}^{k}
\widetilde{U}^{k,n} \stackrel{d}{\longrightarrow} \sum_{k=0}^{\infty
}\frac{a_{k}}{k!} (\sigma B^{H} )^{\diamond k}.
\]
With the choice $a_{n,k}:=1$ and
\[
\sigma_{n} := \frac{\sigma}{1+\mu/n} \rightarrow\sigma \qquad\mbox{as } n \rightarrow\infty,
\]
we observe by (\ref{recursivesystemwickdifferenceforWickpowers}) that $V^{n}_{l} := \widetilde{V}^{n}_{l/n}$
satisfies
\[
V_{l}^{n} = V_{l-1}^{n} + \biggl(\frac{\sigma}{1+\mu/n}
\biggr)V_{l-1}^{n}\diamond_{n} \bigl(B_{l/n}^{H,n} - B_{(l-1)/n}^{H,n} \bigr),\qquad V_{0}^{n} = 1,\qquad l = 1, \ldots, n.
\]
Consider now the piecewise constant function
$ (\widetilde{W}^{n}_{t} )_{t \in[0,1]}$ determined by
$\widetilde{W}^{n}_{t} := W^{n}_{ \lfloor nt \rfloor}$ and
\[
W_{l}^{n} = \biggl(1+\frac{\mu}{n} \biggr)W_{l-1}^{n},\qquad W_{0}^{n} = s_{0},\qquad l = 1,\ldots, n.
\]
By this well-known Euler scheme,
%
\begin{equation}\label{WinL2conv}
(\widetilde{W}^{n}_{t} )_{t\in[0,1]} \longrightarrow s_{0} (\exp(\mu
t) )_{t\in[0,1]}
\end{equation}
in the sup-norm on $[0,1]$.
The product
\begin{eqnarray*}
V_{l}^{n}W_{l}^{n} &=& \biggl(1+\frac{\mu}{n} \biggr)V^{n}_{l-1}W^{n}_{l-1} + \biggl[
\biggl(\frac{\sigma}{1+\mu/n} \biggr)V^{n}_{l-1}\diamond_{n} \bigl(B_{l/n}^{H,n} - B_{(l-1)/n}^{H,n} \bigr) \biggr]
\biggl(1+\frac{\mu}{n}\biggr)W_{l-1}^{n}\\
&=& \biggl(1+\frac{\mu}{n} \biggr)V_{l-1}^{n}W_{l-1}^{n} + \sigma
V_{l-1}^{n}W_{l-1}^{n} \diamond_{n} \bigl(B_{l/n}^{H,n} - B_{(l-1)/n}^{H,n} \bigr),\\
V_{0}^{n}W_{0}^{n} &=& s_{0},\qquad l = 1, \ldots, n,
\end{eqnarray*}
satisfies the Wick difference equation (\ref{discretelinearWickdifferencedrift}) for $S^{n}_{l} = V_{l}^{n}W_{l}^{n}$. The
multiplication by the deterministic function $s_{0}\exp(\mu t)$ is
continuous on the Skorokhod space. Thus, with (\ref{WinL2conv}) and
Billingsley
\cite{Billialt}, Theorem 4.1, we obtain
\[
(\widetilde{S}^{n}_{t} )_{t \in[0,1]} = (\widetilde
{V}^{n}_{t}\widetilde{W}^{n}_{t} )_{t \in[0,1]} \stackrel
{d}{\longrightarrow} s_{0} (\exp(\mu t ) \exp^{\diamond} (\sigma
B^{H}_{t} ) )_{t \in[0,1]}
\]
in the Skorokhod space $D([0,1],\mathbb{R})$. As $s_{0}\exp(\mu t )\exp
^{\diamond} (\sigma B^{H}_{t} )$ solves the SDE (\ref{wicklinearSDEdrift}) (cf.~Mishura \cite{Mishura}, Theorem 3.3.2), the proof is complete.
\end{pf}

\begin{remark}
Theorem \ref{approxlinearSDEwithdrift} holds with additional approximations
$(\sigma_{n},\mu_n) \rightarrow(\sigma,\mu)$.
\end{remark}

\begin{remark}
Theorem \ref{approxlinearSDEwithdrift} was conjectured by Bender
and Elliott \cite{BenderElliott} in their study of the discrete
Wick-fractional Black--Scholes market. They deduced an arbitrage in
this model for sufficiently large $n$. Although the arbitrage or
no-arbitrage property is not preserved by weak convergence, this model
showed that it is even possible to obtain arbitrage in this simple
discrete Wick fractional market model. In a recent work \cite{Valkeila}, Valkeila shows that an alternative approximation to the
exponential of a fractional Brownian motion by a superposition of some
independent renewal reward processes leads to an arbitrage-free and
complete model. We refer to Gaigalas and Kaj \cite{Gaigalas} for a
general limit discussion for these superposition processes.
\end{remark}

\begin{Theo}[(Linear system of SDEs)]\label{examplelinearsystemofsdes}
The piecewise constant interpolation
\[
(\widetilde{X}^{n}_{t},\widetilde{Y}^{n}_{t} )^{\mathrm{T}} := \bigl(X^{n}_{ \lfloor
nt \rfloor},Y^{n}_{ \lfloor nt \rfloor} \bigr)^{\mathrm{T}}
\]
for the solution of the linear system of Wick difference equations
\begin{eqnarray*}
X^{n}_{l} &=& X^{n}_{l-1} + (A_{1}X^{n}_{l-1} + A_{2}Y^{n}_{l-1} )
\diamond_{n} \bigl(B^{H,n}_{l/n} - B^{H,n}_{(l-1)/n} \bigr) ,\qquad
X^{n}_{0} = x_{0},\qquad l = 1, \ldots, n,\\
Y^{n}_{l} &=& Y^{n}_{l-1} + (B_{1}X^{n}_{l-1} + B_{2}Y^{n}_{l-1} )
\diamond_{n} \bigl(B^{H,n}_{l/n} - B^{H,n}_{(l-1)/n} \bigr) ,\qquad
Y^{n}_{0} = y_{0},\qquad l = 1, \ldots, n,
\end{eqnarray*}
converges weakly to the solution $(X,Y)^{\mathrm{T}}$ of the corresponding
linear system of SDEs (\ref{lineare_system})
in the Skorokhod space $D([0,1],\mathbb{R})^2$.
\end{Theo}

\begin{pf}
Analogously to (\ref{linearSDEsystemsolution}), we obtain, by the
recursive system of Wick difference equations for $U^{k,n}$ in
(\ref{recursivesystemwickdifferenceforWickpowers}), the
coefficients for the solution of the systems of difference equations
\[
X_{l}^{n} = \sum_{k =0}^{\infty}\frac{a_{k}}{k!} U^{k,n}_{l} ,\qquad
Y_{l}^{n} = \sum_{k =0}^{\infty}\frac{b_{k}}{k!} U^{k,n}_{l},
\]
recursively by
\[
a_{0} = x_{0} ,\qquad b_{0} = y_{0} ,\qquad a_{k} = A_{1}a_{k-1} + A_{2}b_{k-1}
,\qquad
b_{k} = B_{1}a_{k-1} + B_{2}b_{k-1}.
\]
We define the upper bound
\[
M_{AB} := 2\max\{|A_{1}|, |A_{2}|, |B_{1}|, |B_{2}| \}.
\]
Suppose that $r_{1}, r_{2} \in\mathbb{R}$ are arbitrary. By the linear
system and (\ref{recursivesystemwickdifferenceforWickpowers}),
the sequence of processes
\[
r_{1}\widetilde{X}^{n} + r_{2}\widetilde{Y}^{n} = \sum_{k = 0}^{n}
\biggl(\frac{r_{1}a_{k} + r_{2}b_{k}}{k!} \biggr)\widetilde{U}^{k,n}
\]
fulfils the conditions in Theorem \ref{extendedmaintheorem} with
\[
|r_{1}a_{k} + r_{2}b_{k}| \leq\max\{|x_{0}|,|y_{0}| \}
(|r_{1}|+|r_{2}| )M_{AB}^{k}.
\]
Thus, we obtain the weak convergence
\[
r_{1}\widetilde{X}^{n} + r_{2}\widetilde{Y}^{n} \stackrel
{d}{\longrightarrow} \sum_{k =0}^{\infty} \biggl(\frac{r_{1}a_{k} +
r_{2}b_{k}}{k!} \biggr) (B^{H,n} )^{\diamond k} =r_{1}X + r_{2}Y.
\]
The Cram\'{e}r--Wold device (Billingsley \cite{Billialt}, Theorem 7.7) can now be used to complete the proof.
\end{pf}

\begin{remark}
Theorem \ref{examplelinearsystemofsdes} can be extended to
higher-dimensional linear cases. It also holds for an additional
approximation of the coefficients $A_{n,i} \rightarrow A_{i}$ and
$B_{n,i} \rightarrow B_{i}$ for $n \rightarrow\infty$.
\end{remark}

\begin{example}[(Wick-sine and Wick-cosine)]
The piecewise constant interpolation of
\begin{eqnarray*}
X_{l}^{n} &=& X_{l-1}^{n} + Y_{l-1}^{n}\diamond_{n} \bigl(B^{H,n}_{l/n} - B^{H,n}_{(l-1)/n} \bigr),\qquad X^{n}_{0} = 0,\qquad l = 1, \ldots, n,\\
Y_{l}^{n} &=& Y_{l-1}^{n} - X_{l-1}^{n}\diamond_{n} \bigl(B^{H,n}_{l/n} - B^{H,n}_{(l-1)/n} \bigr),\qquad Y^{n}_{0} = 1,\qquad l = 1, \ldots,
n,
\end{eqnarray*}
converges weakly to the solution of the linear system
\begin{eqnarray*}
\mathrm{d}X_{t} &=& Y_{t}\,\mathrm{d}^{\diamond}B^{H}_{t},\hspace*{10pt}\qquad X_{0} = 0,\\
\mathrm{d}Y_{t} &=& - X_{t}\,\mathrm{d}^{\diamond}B^{H}_{t},\qquad Y_{0} = 1,
\end{eqnarray*}
the process $ (\sin^{\diamond} (B^{H}_{t} ) , \cos^{\diamond}
(B^{H}_{t} ) )^{\mathrm{T}}$. By Theorem \ref{dattheo1},
it can also be approximated by the discrete Wick version functional
$ (\sin^{\diamond_{n}}(B^{H,n}_{t}) , \cos^{\diamond_{n}}(B^{H,n}_{t})
)^{\mathrm{T}}$.
\end{example}

\section{Walsh decompositions and $L^2$-estimates}\label{WickDiffsec}

In this section, we give the Walsh decompositions for the
approximating sequences and obtain some $L^2$-estimates. A key to
the approximation results will be the convergence of the $L^2$-norms
of the discrete Wick powers of $B^{H,n}_{t}$ to the corresponding
$L^2$-norms of the Wick powers of $B^{H}_{t}$.

Recall the Walsh decomposition $B_{t}^{H,n} = \sum^{ \lfloor nt \rfloor
}_{i = 1} b_{t,i}^{n} \xi^{n}_{i}$.
Define
\[
b^{n}_{t,A} := \prod_{i \in A} b_{t,i}^{n},\qquad \Xi^{n}_{A} := \prod_{i \in
A} \xi^{n}_{i},\qquad d_{l,i}^{n} := b_{l/n,i}^{n} - b_{(l-1)/n,i}^{n}
\]
for $l = 1, \ldots, n$. Note that $d_{i,i}^{n} = b_{i/n,i}^{n}$,
$d_{l,i}^{n} = 0$, for $i > l$ and that the increment has the representation
\[
B_{l/n}^{H,n} - B_{(l-1)/n}^{H,n} = \sum_{i = 1}^{l}
d_{l,i}^{n}\xi_{i}^{n}.
\]

Recall the recursive system of Wick difference equations,
%
\begin{equation}\label{dassystemderdiffeq}
U_{l}^{k,n} = U_{l-1}^{k,n} + k U_{l-1}^{k-1,n}\diamond_{n} \bigl(B_{l/n}^{H,n} - B_{(l-1)/n}^{H,n} \bigr) ,\qquad U^{0,n}_{l} = 1
,\qquad U^{k,n}_{0} = 0,
\end{equation}
for $l = 1, \ldots, n$ and $k \in\mathbb{N}$.

\begin{Prop}\label{theowalshdiscretewickpowers}
For all $n,k \in\mathbb{N}$ and $l = 0, \ldots, n$, we have the Walsh
decompositions
%
\begin{eqnarray}\label{walshforU}
\frac{1}{k!}U^{k,n}_{l} &=& \mathop{\sum_{C \subseteq\{1,\ldots, l\}
}}_{ |C| = k} \biggl(\mathop{\sum_{m:C \rightarrow\{1,\ldots, l\}}}_{\mathrm{injective}}\prod_{p \in C} d_{m(p),p}^n \biggr)\Xi_{C}^{n},
\\
\label{walshforbinaryrandomwalk}
\frac{1}{k!} (B^{H,n}_{l/n} )^{\diamond_{n} k} &=& \mathop{\sum_{
 C \subseteq\{1,\ldots, l\}}}_{ |C| = k} b_{l/n,C}^{n}\Xi_{C}^{n},
\\
\label{walshforthedifference}
\frac{1}{k!} (B^{H,n}_{l/n} )^{\diamond_{n} k} - \frac
{1}{k!}U^{k,n}_{l} &=& \mathop{\sum_{C \subseteq\{1,\ldots, l\}}}_{
|C| = k} \biggl(\mathop{\sum_{m:C \rightarrow\{1,\ldots, l\} }}_{\mathrm{not\ injective}}\prod_{p \in C} d_{m(p),p}^n \biggr)\Xi_{C}^{n}.
\end{eqnarray}
\end{Prop}

\begin{pf}
We use the conventions that an empty sum is zero, an empty product
is one and that there exists exactly one map from the empty set to
an arbitrary set. For these reasons, the formulas hold for $k = 0$ or
$l = 0$. We prove (\ref{walshforU}) by induction as follows. For all
$l = 0,
\ldots, n$ and all $k \in\mathbb{N}$, it is obvious that $U^{0,n}_{l}
= 1$
and $U^{k,n}_{0} = 0$, as in formula (\ref{walshforU}). Suppose the
formula is proved for all positive integers less than or equal to a certain
$k$ and all $l = 0, \ldots, n$. Furthermore, for $k+1$, suppose the
formula is proved for all positive integers less than or equal to a certain
$l$. For $k+1$ and $l+1$, we compute, by the difference equation
(\ref{dassystemderdiffeq}) and the induction hypothesis,
%
\begin{eqnarray}\label{inUwalsheq1}
U_{l+1}^{k+1,n} - U_{l}^{k+1,n} &=& (k+1)k! \biggl(\mathop{\sum_{C
\subseteq\{1,\ldots, l\} }}_{ |C| = k} \mathop{\sum_{m:C \rightarrow
\{1,\ldots, l\} }}_{ \mathrm{injective}} \prod_{p \in C} d_{m(p),p}^{n} \Xi
_{C}^n \biggr) \diamond_{n} \sum_{i = 1}^{l+1} d_{l+1,i}^{n}\xi
_{i}^{n}\nonumber\\
&=& (k+1)! \mathop{\mathop{\sum_{C \subseteq\{1,\ldots, l\} }}_{ i \in\{
1,\ldots, l+1\} }}_{ |C| = k , i \notin C}\mathop{\sum_{m:C
\rightarrow\{1,\ldots, l\} }}_{ \mathrm{injective}} \prod_{p \in C}
d_{m(p),p}^{n}d_{l+1,i}^{n} \Xi_{C \cup\{i\}}^n\\
&=& (k+1)! \mathop{\sum_{C' \subseteq\{1,\ldots, l+1\} }}_{ |C'| =
k+1} \mathop{\sum_{m':C' \rightarrow\{1,\ldots, l+1\} }}_{
\mathrm{injective} ,\ \exists q : m(q) = l+1}\prod_{p \in C'} d_{m(p),p}^{n}\Xi
_{C'}^n.\nonumber
\end{eqnarray}
Note that $d_{m,p} = 0$ for all $p-1 \geq m$. Thus, by the induction hypothesis,
%
\begin{equation}\label{inUwalsheq2}
U_{l}^{k+1,n} = (k+1)!\mathop{\sum_{C \subseteq\{1,\ldots, l+1\} }}_{
|C| = k+1} \mathop{\sum_{m:C \rightarrow\{1,\ldots, l+1\} }}_{\mathrm{injective}
,\ \forall q : m(q) < l+1} \prod_{p \in C} d_{m(p),p}^{n}
\Xi_{C}^n.
\end{equation}
Thanks to equations (\ref{inUwalsheq1}) and (\ref{inUwalsheq2}), we obtain (\ref{walshforU}).
In particular, $U_{l}^{k,n} = 0$ for all $k > l$. We now compute the
$k$th Wick power of
$ (B^{H,n}_{l/n} )$ as follows:
\begin{eqnarray*}
\Biggl(\sum ^{n}_{i = 1} b^{n}_{l/n,i} \xi^{n}_{i} \Biggr)^{\diamond_{n}
k} &=& \mathop{\sum ^{l}_{i_{1}, i_{2},\ldots, i_{k} = 1 }}_{
\mathrm{pairwise\ distinct}} \Biggl(\prod^{k}_{j = 1} b^{n}_{l/n,i_{j}}
\prod^{k}_{j = 1}\xi^{n}_{i_{j}} \Biggr)\\
&=& \mathop{\sum _{C \subseteq\{1,\ldots, l\} }}_{|C| = k} k! \Biggl(\prod
_{i \in C} b^{n}_{l/n,i} \prod_{i \in C} \xi^{n}_{i} \Biggr) =\mathop{\sum
_{C \subseteq\{1,\ldots, l\}}}_{ |C| = k} k! b^{n}_{l/n,C} \Xi^{n}_{C}.
\end{eqnarray*}
In particular, $ (B^{H,n}_{l/n} )^{\diamond_{n} k} = 0$ for all
$k > l$. This yields (\ref{walshforbinaryrandomwalk}).

The telescoping sum yields
\[
\sum_{m(p)=1}^{l} d_{m(p),p}^{n} = \sum_{m(p)=p}^{l} d_{m(p),p}^{n} =
\sum_{m(p) = p}^{l} \bigl(b_{m(p)/n,p}^{n} - b_{(m(p)-1)/n,p}^{n} \bigr) =b_{l/n,p}^{n}
\]
and thus we get
%
\begin{equation}\label{X_teleskopsumme}
\sum_{m:C \rightarrow\{1,\ldots, l\}} \prod_{p \in C}
d_{m(p),p}^{n} =\prod_{p \in C} \Biggl(\sum_{m(p) = 1}^{l} d_{m(p),p}^{n} \Biggr)
=\prod_{p \in C} b_{l/n,p}^{n} = b_{l/n,C}^{n}.
\end{equation}
Equation (\ref{walshforthedifference}) is, thus, implied by (\ref{walshforU}) and (\ref{walshforbinaryrandomwalk}).
\end{pf}

In the following propositions, we obtain some elementary estimates for
the $L^2$-norm of discrete Wick powers of $B^{H,n}_{t}$.

\begin{Prop}\label{dasheiligewort}
For all $t \in[0,1]$ and $i \in\{1, \ldots, \lfloor nt \rfloor\}$,
\[
b_{t,i}^{n} \leq2c_{H} n^{-(1-H)}.
\]
\end{Prop}

\begin{pf} We estimate
\begin{eqnarray*}
b_{t,i}^{n} &=& n^{1/2}c_{H}\biggl(H-\frac{1}{2}\biggr)\int_{
(i-1)/n}^{i/n} s^{1/2-H} \int_{s}^{\lfloor nt
\rfloor/n}u^{H-1/2} (u-s )^{H-3/2}\,\mathrm{d}u \,\mathrm{d}s\\
&\leq& n^{1/2}c_{H}\biggl(H-\frac{1}{2}\biggr)\int_{(i-1)/n}^{i/n}
s^{1/2-H} \biggl(\frac{ \lfloor nt \rfloor}{n} \biggr)^{H-
1/2}\frac{1}{H-1/2} \biggl(\frac{ \lfloor nt \rfloor}{n}-s
\biggr)^{H-1/2}\, \mathrm{d}s\\
&\leq& n^{1/2}c_{H} \frac{1}{3/2-H}\biggl ( \biggl(\frac{i}{n}
\biggr)^{3/2-H} - \biggl(\frac{i-1}{n} \biggr)^{3/2-H} \biggr) \biggl(\frac{ \lfloor
nt \rfloor}{n} \biggr)^{2(H-1/2)}.
\end{eqnarray*}
Since $t \leq1$, $\frac{1}{3/2-H} \leq2$ and
$|x|^{3/2-H} - |y|^{3/2-H} \leq
|x-y|^{3/2-H}$, the assertion follows.
\end{pf}

\begin{remark}\label{remarkonboundedness}
Observe that
%
\begin{equation}\label{covarianceforbinaryramdomwalk}
\mathbf{E} [B^{H,n}_{t} B^{H,n}_{s} ] = \mathbf{E} \Biggl[\sum_{i_{1}, i_{2}
= 1}^{ \lfloor nt \rfloor} b_{t,i_{1}}^{n} b_{s,i_{2}}^{n}\xi
_{i_{1}}^{n}\xi_{i_{2}}^{n} \Biggr] = \sum^{ \lfloor nt \rfloor}_{i = 1}
(b_{t,i}^{n} b_{s,i}^{n} ).
\end{equation}
By Nieminen \cite{Nieminen}, we thus get, for any $s,t \in[0,1]$, the
following convergence:
\begin{eqnarray}\label{nieminenconvergencefornorms}
\mathbf{E} [B^{H,n}_{t} B^{H,n}_{s} ] &=& \sum^{ \lfloor nt \rfloor}_{i =
1} n \int^{i/n}_{(i-1)/n}z \biggl(\frac{ \lfloor nt \rfloor
}{n}, u \biggr)\,\mathrm{d}u \int^{i/n}_{(i-1)/n}z \biggl(\frac{ \lfloor ns
\rfloor}{n}, u \biggr)\,\mathrm{d}u \nonumber\\[-8pt]\\[-8pt]
&\longrightarrow& \int_{0}^{1}z(t,u)z(s,u)\,\mathrm{d}u =\mathbf{E} [B^{H}_{t}
B^{H}_{s} ].\nonumber
\end{eqnarray}
Moreover, we have, by the Cauchy--Schwarz inequality, the upper bound
%
\begin{eqnarray}\label{zurboundednessnorm}
\mathbf{E} [ (B^{H,n}_{t} - B^{H,n}_{s} )^2 ] &=&\sum^{ \lfloor nt
\rfloor}_{i=1} \biggl(\sqrt{n} \int^{i/n}_{(i-1)/n} \biggl(z\biggl(\frac{
\lfloor nt \rfloor}{n}, u\biggr)\,\mathrm{d}u - z\biggl(\frac{ \lfloor ns \rfloor}{n}, u\biggr) \biggr)\,\mathrm{d}u
\biggr)^2\nonumber\\
&\leq& \sum^{ \lfloor nt \rfloor}_{i=1} \int^{i/n}_{
(i-1)/n} \biggl(z\biggl(\frac{ \lfloor nt \rfloor}{n}, u\biggr) - z\biggl(\frac{ \lfloor ns
\rfloor}{n}, u\biggr) \biggr)^2\,\mathrm{d}u\\
&=& \biggl|\frac{ \lfloor nt \rfloor}{n} - \frac{ \lfloor ns \rfloor}{n}
\biggr|^{2H}.\nonumber
\end{eqnarray}
\end{remark}

\begin{Prop}\label{walshpropevoll}
For all $t \geq s$ in $[0,1]$ and all $N \in\mathbb{N}$ such that $
\lfloor ns \rfloor\geq N$, we have
%
\begin{eqnarray}\label{L2normestimateofthedifference}
0 &\leq& \mathbf{E} [(B^{H,n}_{t})^2 ]^{N} + \mathbf{E}
[(B^{H,n}_{s})^2 ]^{N} - 2\mathbf{E} [B^{H,n}_{t} B^{H,n}_{s}
]^{N}\nonumber
\\
&&{}- \frac{1}{N!}\mathbf{E} \bigl[ \bigl( (B^{H,n}_{t} )^{\diamond_{n} N} -
(B^{H,n}_{s} )^{\diamond_{n} N} \bigr)^2 \bigr] \\
\nonumber\\&\leq& 2c_{H}^{2} N^2 t^{2H (N-1)} n^{-(2-2H)}. \nonumber
\end{eqnarray}
In particular,
\[
\lim_{n\rightarrow
\infty}\mathbf{E} \bigl[ \bigl( (B^{H,n}_{t} )^{\diamond_{n}
N} - (B^{H,n}_{s} )^{\diamond_{n} N} \bigr)^2 \bigr]
= \mathbf{E} \bigl[ \bigl( (B^{H}_{t} )^{\diamond N} - (B^{H}_{s} )^{\diamond
N} \bigr)^2 \bigr].
\]
\end{Prop}

\begin{pf}
As $N \leq \lfloor ns \rfloor$, we get, making use of
Proposition \ref{theowalshdiscretewickpowers}, (\ref{L2norminWalshterms}) and
(\ref{covarianceforbinaryramdomwalk}) in Remark \ref{remarkonboundedness},
\begin{eqnarray}\label{zurnormestimate}
&&\frac{1}{N!}\mathbf{E} [(B^{H,n}_{t})^{\diamond_{n}
N}(B^{H,n}_{s})^{\diamond_{n} N} ]\nonumber\\
&&\qquad= \frac{1}{N!}\mathbf{E} \biggl[ \biggl(N!\mathop{\sum_{C \subseteq \{1, \ldots
, \lfloor nt \rfloor\} }}_{ |C| = N} b_{t,C}^{n}\Xi_{C}^{n} \biggr) \biggl(N!\mathop{\sum
_{C \subseteq \{1, \ldots, \lfloor ns \rfloor\} }}_{ |C| =
N} b_{s,C}^{n}\Xi_{C}^{n} \biggr) \biggr]\nonumber\\
&&\qquad= N!\mathop{\sum_{C \subseteq \{1, \ldots, \lfloor nt \rfloor\}
}}_{|C| = N} b_{t,C}^{n}b_{s,C}^{n}\\
&&\qquad= \sum_{i_{1}, \ldots, i_{N} = 1}^{
\lfloor nt \rfloor} \prod_{j = 1}^{N} (b_{t,i_{j}}^{n} b_{s,i_{j}}^{n}
) - \mathop{\sum_{i_{1}, \ldots, i_{N} = 1}}_{ \exists k,l \dvtx
i_{k} = i_{l}}^{ \lfloor nt \rfloor} \prod_{j = 1}^{N}
(b_{t,i_{j}}^{n} b_{s,i_{j}}^{n} )\nonumber\\
&&\qquad= \mathbf{E} [B^{H,n}_{t} B^{H,n}_{s} ]^{N} - \mathop{\sum_{i_{1},
\ldots, i_{N} = 1}}_{ \exists k,l \dvtx i_{k} = i_{l}}^{ \lfloor nt
\rfloor} \prod_{j = 1}^{N} (b_{t,i_{j}}^{n} b_{s,i_{j}}^{n} ).\nonumber
\end{eqnarray}
Thus, we have
\begin{eqnarray}\label{estimatefornormconv}
&&\mathbf{E} [(B^{H,n}_{t})^2 ]^{N} + \mathbf{E} [(B^{H,n}_{s})^2 ]^{N}
- 2\mathbf{E} [B^{H,n}_{t} B^{H,n}_{s} ]^{N}\nonumber\\
&&\qquad{}- \frac{1}{N!}\mathbf{E} \bigl[ \bigl( (B^{H,n}_{t} )^{\diamond_{n} N} -
(B^{H,n}_{s} )^{\diamond_{n} N} \bigr)^2
\bigr]\nonumber\\[-8pt]\\[-8pt]
&&\quad= \mathop{\sum_{i_{1}, \ldots, i_{N} = 1}}_{ \exists k,l \dvtx i_{k}
= i_{l}}^{ \lfloor nt \rfloor} \Biggl(
\prod_{j = 1}^{N} (b_{t,i_{j}}^{n} )^2 + \prod_{j = 1}^{N}
(b_{s,i_{j}}^{n} )^2 - 2\prod_{j = 1}^{N} (b_{t,i_{j}}^{n}
b_{s,i_{j}}^{n} ) \Biggr)\nonumber\\
&&\quad= \mathop{\sum_{i_{1}, \ldots, i_{N} = 1}}_{
\exists k,l \dvtx i_{k} = i_{l}}^{ \lfloor nt \rfloor}
\Biggl(\prod_{j = 1}^{N} (b_{t,i_{j}}^{n} ) -
\prod_{j = 1}^{N} (b_{s,i_{j}}^{n} ) \Biggr)^2 \geq 0.\nonumber
\end{eqnarray}
Hence, the left-hand side of the inequality in (\ref{L2normestimateofthedifference}) follows.
By Proposition \ref{dasheiligewort}, (\ref{covarianceforbinaryramdomwalk}) and
(\ref{zurboundednessnorm}) in Remark \ref{remarkonboundedness}, as well as $ |\frac{ \lfloor nt \rfloor}{n} | \leq t$, we obtain
%
\begin{eqnarray}\label{restefornorm}
&&\mathop{\sum_{i_{1}, \ldots, i_{N} = 1}}_{ \exists k,l \dvtx i_{k}
= i_{l}}^{ \lfloor nt \rfloor}
\Biggl(\prod_{j = 1}^{N} (b_{t,i_{j}}^{n} ) - \prod_{j = 1}^{N}
(b_{s,i_{j}}^{n} ) \Biggr)^2\nonumber\\
&&\quad\leq \mathop{\sum_{i_{1}, \ldots, i_{N} = 1}}_{ \exists k,l \dvtx
i_{k} = i_{l}}^{ \lfloor nt \rfloor}
\Biggl(\prod_{j = 1}^{N} (b_{t,i_{j}}^{n} ) \Biggr)^2 \leq \pmatrix{ N \cr 2}
\Bigl(\max_{i}(b_{t,i}^{n})^2 \Bigr) \sum_{i_{1}, \ldots, i_{N-1} = 1}^{
\lfloor nt \rfloor} \Biggl(\prod_{j = 1}^{N-1} (b_{t,i_{j}}^{n} ) \Biggr)^2
\\
&&\quad\leq 2c_{H}^{2} N^2\mathbf{E} [ (B^{H,n}_{t} )^2 ]^{N-1} n^{-(2-2H)}
\leq2c_{H}^{2}N^2 t^{2H (N-1)} n^{-(2-2H)} \rightarrow 0\nonumber
\end{eqnarray}
for $n \rightarrow\infty$. The representation of Wick powers of
$B^{H}_{t}$ by Hermite polynomials, as in (\ref{propwickpowersarehermitepolynomials}), their orthonormality (cf.~Kuo \cite{Kuo},
page 355) and the polarization identity collectively yield $\mathbf{E}
[ (B^{H}_{t} )^{\diamond N} (B^{H}_{s} )^{\diamond N} ] = N!\mathbf{E}
[ (B^{H}_{t} ) (B^{H}_{s} ) ]^{N}$
(cf.~also \cite{Nualart}, Lemma 1.1.1).
Thus, we have, by (\ref{estimatefornormconv}),
\begin{eqnarray*}
&&\mathbf{E} \bigl[ \bigl( (B^{H,n}_{t} )^{\diamond_{n} N} - (B^{H,n}_{s}
)^{\diamond_{n} N} \bigr)^2 \bigr] - \mathbf{E} \bigl[ \bigl( (B^{H}_{t} )^{\diamond N} -
(B^{H}_{s} )^{\diamond N} \bigr)^2 \bigr]\nonumber\\
&&\quad= N! \bigl(\mathbf{E} [(B^{H,n}_{t})^2 ]^{N} - \mathbf{E} [(B^{H}_{t})^2
]^{N} + \mathbf{E} [(B^{H,n}_{s})^2 ]^{N} - \mathbf{E} [(B^{H}_{s})^2
]^{N} \nonumber\\
&&\qquad\hphantom{N! \bigl(}{} - 2\mathbf{E} [B^{H,n}_{t} B^{H,n}_{s} ]^{N} +2 \mathbf{E} [B^{H}_{t}
B^{H}_{s} ]^{N} \bigr) \\
&&\qquad{} - N! \mathop{\sum_{i_{1}, \ldots, i_{N} = 1}}_{ \exists k,l \dvtx
i_{k} = i_{l}}^{ \lfloor nt \rfloor}
\Biggl(\prod_{j = 1}^{N} (b_{t,i_{j}}^{n} ) - \prod_{j = 1}^{N}
(b_{s,i_{j}}^{n} ) \Biggr)^2.
\end{eqnarray*}
Applying the convergences (\ref{restefornorm}) and (\ref{nieminenconvergencefornorms}) yields
\[
\mathbf{E} \bigl[ \bigl( (B^{H,n}_{t} )^{\diamond_{n} N} - (B^{H,n}_{s}
)^{\diamond_{n} N} \bigr)^2 \bigr] - \mathbf{E} \bigl[ \bigl( (B^{H}_{t} )^{\diamond N} -
(B^{H}_{s} )^{\diamond N} \bigr)^2 \bigr] \rightarrow 0.
\]
\upqed
\end{pf}

\begin{remark}\label{remarkonnormofwickpowers}
In particular, we obtain, by (\ref{zurnormestimate}), (\ref{zurboundednessnorm}) and $ |\frac{ \lfloor nt \rfloor}{n} | \leq t$,
\begin{eqnarray}\label{boundedness_norm}
\mathop{\sum_{C \subseteq\{1,\ldots, \lfloor nt \rfloor\} }}_{ |C| =
N} (b_{t,C}^{n} )^2 &=& \biggl(\frac{1}{N!} \biggr)^2 \mathbf{E} [ ( (B^{H,n}_{t}
)^{\diamond_{n} N} )^2 ] \nonumber\\[-8pt]\\[-8pt]
&\leq& \frac{1}{N!}\mathbf{E} [ (B^{H,n}_{t} )^2 ]^{N} \leq\frac
{1}{N!} \biggl|\frac{ \lfloor nt \rfloor}{n} \biggr|^{2H N} \leq\frac{1}{N!} t^{2H
N}.\nonumber
\end{eqnarray}
\end{remark}

The next proposition estimates the difference between the approximating
sequences in Theorems \ref{dattheo1} and \ref{extendedmaintheorem}.

\begin{Prop}\label{propzumdiffconvrate}
Under the assumptions of Theorem \ref{dattheo1}, there exists a
constant $K > 0$ such that for all $t \in[0,1]$, $n \geq1$ and $k \in
\mathbb{N}$,
%
\begin{equation}\label{diediffconvrateimextendedmaintheorem}
\mathbf{E} \Biggl[ \Biggl|\sum_{k=0}^{n}\frac{a_{n,k}}{k!} (B_{t}^{H,n} )^{\diamond
_{n} k} - \sum_{k=0}^{n}\frac{a_{n,k}}{k!} \widetilde{U}^{k,n}_{t} \Biggr|^2
\Biggr] \leq K n^{1-2H}
\end{equation}
for the approximating processes in Theorems \ref{dattheo1} and \ref{extendedmaintheorem}.
\end{Prop}

\begin{pf}
Recall that
$d_{r,i}^{n} =b_{r/n,i}^{n} - b_{(r-1)/n,i}^{n} =\sqrt
{n}\int_{(i-1)/n}^{i/n} (z(\frac{r}{n},s) - z(\frac
{r-1}{n},s) )\,\mathrm{d}s$.
By (\ref{zurboundednessnorm}) in Remark \ref{remarkonboundedness},
we obtain
\[
(d_{r,i}^{n} )^2 \leq\sum_{i = 1}^{r} (d_{r,i}^{n} )^2 \leq \biggl|\frac
{r}{n} - \frac{r-1}{n} \biggr|^{2H} =n^{-2H}.
\]
Thus, we have $d_{r,i}^{n} \leq n^{-H}$ for all $i,n,r \geq1$. Hence,
we obtain, as the sum in (\ref{X_teleskopsumme}) telescopes, for $|C|
\geq2$,
%
\begin{eqnarray}\label{anestimateforZidinger}
\mathop{\sum_{m:C \rightarrow\{1,\ldots, \lfloor nt \rfloor\} }}_{
\mathrm{not\ injective}}\prod_{l \in C} d_{m(l),l}^n &=&\mathop{\sum_{
m:C \rightarrow\{1,\ldots, \lfloor nt \rfloor\} }}_{ \exists u,v \in C
\dvtx m(u) = m(v) } \prod_{l \in C} d_{m(l),l}^n \nonumber\\
&=& \sum_{u \in C} \sum_{m:C \setminus\{u\} \rightarrow\{1,\ldots,
\lfloor nt \rfloor\}} \biggl(\prod_{l \in C \setminus\{u\}} d_{m(l),l}^n
\biggr)\sum_{v \in C\setminus\{u\}} d_{m(v),u}^n\\
&\leq&\max_{i,r} d_{r,i}^{n} (|C|-1 ) \mathop{\sum_{C' \subset C }}_{
|C'| = |C|-1}\sum_{m:C' \rightarrow\{1,\ldots, \lfloor nt \rfloor\} }
\prod_{l \in C'} d_{m(l),l}^n\nonumber\\
&\leq& n^{-H} (|C|-1) \mathop{\sum_{C' \subset C }}_{ |C'| = |C|-1 }
b_{t,C'}^{n}.\nonumber
\end{eqnarray}
By (\ref{walshforthedifference}), (\ref{anestimateforZidinger}), (\ref{boundedness_norm}) and since $ ( \lfloor nt \rfloor
-(k-1) ) \leq n$, we obtain, for $k \geq1$,
\begin{eqnarray*}
&&\mathbf{E} \biggl[ \biggl|\frac{1}{k!} \bigl( (B_{t}^{H,n} )^{\diamond_{n} k} -
\widetilde{U}^{k,n}_{t} \bigr) \biggr|^2 \biggr]\\
&&\quad\leq \mathop{\sum_{C \subseteq\{1,\ldots, \lfloor nt
\rfloor\}}}_{|C| = k} \biggl( n^{-H} (k-1) \mathop{\sum_{C' \subset C }}_{ |C'| = |C|-1}
b_{t,C'}^{n} \biggr)^2\\
&&\quad\leq n^{-2H} (k-1)^2 \mathop{\sum_{C \subseteq\{1,\ldots, \lfloor
nt \rfloor\} }}_{ |C| = k} (k-1) \mathop{\sum_{C' \subset C }}_{ |C'| =
|C|-1} (b_{t,C'}^{n} )^2\\
&&\quad\leq n^{-2H} (k-1)^3 \bigl( \lfloor nt \rfloor-(k-1) \bigr) \mathop{\sum_{C'
\subseteq\{1,\ldots, \lfloor nt \rfloor\} }}_{ |C'| = k-1}
(b_{t,C'}^{n} )^2 \leq\frac{(k-1)^3}{(k-1)!}t^{2H (k-1)} n^{1-2H}.
\end{eqnarray*}
Since $a_{n,k} \leq C^{k}$ in Theorems \ref{dattheo1} and \ref{extendedmaintheorem}, and $ ( (B_{t}^{H,n} )^{\diamond_{n} k} -
\widetilde{U}^{k,n}_{t} )$ are zero for $k = 0,1$ and orthogonal for
different $k$ by Proposition \ref{theowalshdiscretewickpowers} and
(\ref{L2norminWalshterms}), we get
\begin{eqnarray*}
&&\mathbf{E} \Biggl[ \Biggl|\sum_{k=0}^{n}\frac{a_{n,k}}{k!} (B_{t}^{H,n} )^{\diamond
_{n} k} - \sum_{k=0}^{n}\frac{a_{n,k}}{k!} \widetilde{U}^{k,n}_{t} \Biggr|^2
\Biggr]\\
&&\quad= \sum_{k = 2}^{ \lfloor nt \rfloor} \mathbf{E} \biggl[ \biggl|\frac{a_{n,k}}{k!} \bigl(
(B_{t}^{H,n} )^{\diamond_{n} k} - \widetilde{U}^{k,n}_{t} \bigr) \biggr|^2 \biggr] \leq
\Biggl(\sum_{k = 2}^{\infty} C^{2k} \frac{(k-1)^3}{(k-1)!} t^{2H(k-1)} \Biggr) n^{1-2H}.
\end{eqnarray*}
As the series on the right-hand side converges uniformly in $t \in
[0,1]$, the assertion follows.
\end{pf}

\section{Convergence of the finite-dimensional distributions}\label{secconvfd}

We first prove that Theorems \ref{dattheo1} and \ref{extendedmaintheorem} hold with weak convergence replaced by convergence
of the finite-dimensional distributions. To this end, we first
approximate the Wick powers of $B^{H}_{t}$ by induction. We then
combine these convergence results to approximate the Wick analytic
functionals $F^{\diamond}(B^{H}_{t}) = \sum_{k = 0}^{\infty} \frac
{a_{k}}{k!} (B^{H}_{t} )^{\diamond k}$. Finally, we conclude that
convergence in finite dimensions holds in Theorem \ref{extendedmaintheorem}.

We observed in Section \ref{Wickexpsection} that
$(B^{H}_{t})^{\diamond N} = h^{N}_{|t|^{2H}}(B^{H}_{t})$
and that the Hermite recursion formula
%
\begin{equation}\label{hermiterekursionfBB}
(B^{H}_{t})^{\diamond(N+1)}=(B^{H}_{t})(B^{H}_{t})^{\diamond N} -
|t|^{2H} N (B^{H}_{t})^{\diamond(N-1)}
\end{equation}
holds. For the discrete Wick powers of the discrete variables, we now
obtain a discrete variant of (\ref{hermiterekursionfBB}).

\begin{Prop}[(Discrete Hermite recursion)]\label{lemma_hermite}
For all $N \geq1$ and $t \in[0,1]$,
%
\begin{equation}\label{hermite_discrete}
(B^{H,n}_{t})^{\diamond_{n} (N+1)} = B^{H,n}_{t} (B^{H,n}_{t})^{\diamond
_{n} N} - N \mathbf{E} [(B^{H,n}_{t})^2 ] (B^{H,n}_{t})^{\diamond_{n}
(N-1)} +R(B^{H,n}_{t},N),
\end{equation}
with remainder
%
\begin{equation}\label{R_rest}
R(B^{H,n}_{t},N) = N! \mathop{\sum_{C \subseteq\{1,\ldots \lfloor
nt \rfloor\} }}_{ |C| = N-1} b_{t,C}^{n} \Xi_{C}^{n} \sum_{i \in C}
(b_{t,i}^{n} )^2
\end{equation}
and
%
\begin{equation}\label{restinL2}
\mathbf{E} [ (R(B^{H,n}_{t},N) )^2 ] \leq16 c_{H}^{4} N!N^3 n^{-(4-4H)}.
\end{equation}
\end{Prop}

In particular, we will use the fact that the discrete Hermite recursion
(\ref{hermite_discrete}) converges weakly to Hermite recursion (\ref{hermiterekursionfBB}) for $n \rightarrow\infty$.

\begin{pf}
By Proposition \ref{theowalshdiscretewickpowers}, we get
\begin{eqnarray}\label{hermite_reste}
(B^{H,n}_{t})^{\diamond_{n} (N+1)} &=& B^{H,n}_{t} \diamond_{n}
(B^{H,n}_{t})^{\diamond_{n} N} = \Biggl(\sum^{ \lfloor nt \rfloor}_{i = 1}
b_{t,i}^{n} \xi^{n}_{i} \Biggr) \diamond_{n} \biggl(\mathop{\sum_{A \subseteq\{
1,\ldots, \lfloor nt \rfloor\} }}_{ |A| = N} N! b_{t,A}^{n}
\Xi^{n}_{A}\biggr)\nonumber\\[-8pt]\\[-8pt]
&=& B^{H,n}_{t} (B^{H,n}_{t})^{\diamond_{n} N} - \mathop{\sum_{A
\subseteq\{1,\ldots, \lfloor nt \rfloor\} }}_{ |A| = N} \sum_{i \in A}
N! b_{t,i}^{n} b_{t,A}^{n} \Xi^{n}_{A}\xi^{n}_{i}.\nonumber
\end{eqnarray}
For the second term in the last line in equation (\ref{hermite_reste}), by (\ref{covarianceforbinaryramdomwalk})
in Remark \ref{remarkonboundedness} and Proposition \ref{theowalshdiscretewickpowers}, we obtain
\begin{eqnarray*}
&&\mathop{\sum_{A \subseteq\{1,\ldots, \lfloor nt \rfloor\} }}_{ |A| =
N} \sum_{i \in A} N! b_{t,i}^{n} b_{t,A}^{n} \Xi^{n}_{A}\xi^{n}_{i}\\
&&\quad = N!\mathop{\sum_{A \subseteq\{1,\ldots, \lfloor nt \rfloor\}
}}_{|A| = N} \sum_{i \in A} b_{t,A \setminus\{i\}}^{n} \Xi^{n}_{A
\setminus\{i\}} (b_{t,i}^{n}\xi^{n}_{i} )^2 = N! \mathop{\sum_{C
\subseteq\{1,\ldots, \lfloor nt \rfloor\} }}_{ |C| = N-1} b_{t,C}^{n}
\Xi_{C}^{n} \sum_{i \notin C} (b_{t,i}^{n} )^2\\
&&\quad = N (N-1)! \mathop{\sum_{C \subseteq\{1,\ldots, \lfloor nt \rfloor
\} }}_{ |C| = N-1} b_{t,C}^{n} \Xi_{C}^{n} \Biggl(\sum_{i = 1} ^{ \lfloor nt
\rfloor} (b_{t,i}^{n} )^2 - \sum_{i \in C} (b_{t,i}^{n} )^2 \Biggr)\\
&&\quad= N (B^{H,n}_{t})^{\diamond_{n} (N-1)} \mathbf{E} [(B^{H,n}_{t})^2 ] -
N! \mathop{\sum_{C \subseteq\{1,\ldots, \lfloor nt \rfloor\} }}_{
|C| = N-1} b_{t,C}^{n} \Xi_{C}^{n} \sum_{i \in C} (b_{t,i}^{n} )^2,
\end{eqnarray*}
which yields (\ref{hermite_discrete}) and (\ref{R_rest}).
Thus, thanks to Proposition \ref{dasheiligewort}, Remark \ref{remarkonnormofwickpowers} and $t \leq1$, we obtain
\begin{eqnarray*}
\mathbf{E} [ (R(B^{H,n}_{t},N) )^2 ] &=&(N!)^2 \mathop{\sum_{C
\subseteq\{1,\ldots, \lfloor nt \rfloor\} }}_{ |C| = N-1}
(b_{t,C}^{n})^2 \biggl(\sum_{i \in C} (b_{t,i}^{n} )^2 \biggr)^2\\
&\leq& (N!)^2 \frac{1}{(N-1)!}t^{2H (N-1)} \bigl((N-1)4c_{H}^2 n^{-(2-2H)}
\bigr)^2\\
& \leq&16c_{H}^{4} N! N^3 n^{-(4-4H)}.
\end{eqnarray*}
\upqed
\end{pf}

\begin{Theo}\label{dasTheoalleWickschwach}
For all $N \in\mathbb{N}$,
%
\begin{equation}\label{wickpowerfd}
(1,B^{H,n},\ldots,(B^{H,n})^{\diamond_{n} N} ) \stackrel
{fd}{\longrightarrow} (1,B^{H},\ldots,(B^{H})^{\diamond N} ).
\end{equation}
\end{Theo}

\begin{pf}
The proof proceeds by induction on $N$. By Sottinen's
approximation, $ (1, B^{H,n} ) \stackrel{fd}{\longrightarrow} (1, B^{H}
)$. Suppose that equation (\ref{wickpowerfd}) is proved for some $N
\geq1$. Assume that $k \in\mathbb{N}$ and $r_{i}^{j} \in\mathbb{R}$
for $j = 0, \ldots, N+1$ , $i = 1 , \ldots, k$ and $t_{1}, t_{2},
\ldots, t_{k} \in[0,1]$ are chosen arbitrarily. By the pointwise
convergence $\mathbf{E} [(B^{H,n}_{t})^2 ] \rightarrow|t|^{2H}$ and
the generalized continuous mapping theorem (Billingsley \cite{Billialt}, Theorem~5.5), the induction hypothesis implies
that
\begin{eqnarray*}
&&\sum_{l = 0}^{N} \Biggl(\sum_{j = 1}^{k} r^{l}_{j} (B_{t_{j}}^{H,n})^{\diamond
_{n} l} \Biggr) + \sum_{j = 1}^{k} r^{N+1}_{j} \bigl(B^{H,n}_{t_{j}}
(B_{t_{j}}^{H,n})^{\diamond_{n} N} - N \mathbf{E} [(B^{H,n}_{t_{j}})^2
] (B^{H,n}_{t_{j}})^{\diamond_{n} (N-1)} \bigr)\\
&&\quad\stackrel{d}{\longrightarrow} \sum_{l = 0}^{N} \Biggl(\sum_{j = 1}^{k}
r^{l}_{j}(B_{t_{j}}^{H})^{\diamond l} \Biggr) + \sum_{j = 1}^{k} r^{N+1}_{j}
\bigl(B^{H}_{t_{j}} (B_{t_{j}}^{H})^{\diamond N} - N
|t_{j}|^{2H}(B^{H}_{t_{j}})^{\diamond(N-1)} \bigr).
\end{eqnarray*}
Since $H > \frac{1}{2}$, (\ref{restinL2}) yields $R(B^{H,n}_{t},N)
\longrightarrow0$ in $L^{2}(\Omega,P)$.
Thus, by Slutsky's theorem \cite{Billialt}, Theorem 4.1, and the
Hermite recursions (\ref{hermiterekursionfBB}) and (\ref{hermite_discrete}), we obtain
\[
\sum_{l = 0}^{N+1} \Biggl(\sum_{j = 1}^{k}
r^{l}_{j}(B_{t_{j}}^{H,n})^{\diamond_{n} l} \Biggr)
\stackrel{d}{\longrightarrow} \sum_{l = 0}^{N+1} \Biggl(\sum_{j = 1}^{k}
r^{l}_{j}(B_{t_{j}}^{H})^{\diamond l} \Biggr).
\]
By the Cram\'{e}r--Wold device (Billingsley \cite{Billialt}, Theorem 7.7), we have
\[
(1,B^{H,n},\ldots,(B^{H,n})^{\diamond_{n} N+1} ) \stackrel
{fd}{\longrightarrow} (1,B^{H},\ldots,(B^{H})^{\diamond N+1} )
\]
and the induction is complete.
\end{pf}

\begin{Prop}\label{fdmaintheorem}
In the context of Theorem \ref{dattheo1}, convergence holds in
finite-dimensional distributions.
\end{Prop}

\begin{pf}
By Billingsley \cite{Billialt}, Theorem 4.2, it suffices to show that
the following three conditions hold:
%
\begin{eqnarray}\label{geomwconvstraffganz1}
&\displaystyle\forall m \in\mathbb{N},\qquad \sum_{k = 0}^{m} \frac
{a_{n,k}}{k!}(B^{H,n}_{t})^{\diamond_{n} k} \stackrel
{fd}{\longrightarrow} \sum_{k = 0}^{m} \frac
{a_{k}}{k!}(B^{H}_{t})^{\diamond k} \qquad\mbox{as } n
\rightarrow\infty;&
\\
\label{geomwconvstraffganz2}
&\displaystyle\forall t \in[0,1],\qquad \lim_{m \rightarrow\infty} \limsup_{n \rightarrow
\infty} \mathbf{E} \Biggl[\Biggl|\sum_{k = 0}^{n}\frac{a_{n,k}}{k!} (B_{t}^{H,n}
)^{\diamond_{n} k} - \sum_{k = 0}^{m} \frac
{a_{n,k}}{k!}(B^{H,n}_{t})^{\diamond_{n} k}\Biggr| \wedge1 \Biggr] =
0;&\qquad
\\
\label{geomwconvstraffganz3}
&\displaystyle\sum_{k = 0}^{m} \frac{a_{k}}{k!} (B^{H}_{t})^{\diamond k} \stackrel
{fd}{\longrightarrow} \sum_{k = 0}^{\infty} \frac{a_{k}}{k!}
(B^{H}_{t})^{\diamond k} \qquad\mbox{as } m \rightarrow\infty.&
\end{eqnarray}
Condition (\ref{geomwconvstraffganz1}) follows directly from
Theorem \ref{dasTheoalleWickschwach} and the generalized continuous
mapping theorem (\cite{Billialt}, Theorem 5.5). For the second
condition, we compute
\begin{eqnarray*}
&&\mathbf{E} \Biggl[ \Biggl(\sum_{k = 0}^{n} \frac{a_{n,k}}{k!}
(B^{H,n}_{t})^{\diamond k} - \sum_{k = 0}^{m} \frac
{a_{n,k}}{k!}(B^{H,n}_{t})^{\diamond_{n} k} \Biggr)^2 \Biggr] \\
&&\quad= \mathbf{E} \Biggl[ \Biggl(\sum
_{k = m+1}^{n} \frac{a_{n,k}}{k!}(B_{t}^{H,n})^{\diamond_{n} k} \Biggr)^{2}
\Biggr]\\
&&\quad= \sum_{k = m+1}^{n} \biggl(\frac{a_{n,k}}{k!} \biggr)^2 \mathbf{E} [
((B_{t}^{H,n})^{\diamond_{n} k} )^{2} ] \leq \sum_{k = m+1}^{n} \biggl(\frac
{C^{k}}{k!} \biggr)^2 k! t^{2Hk},
\end{eqnarray*}
applying the estimate of Remark \ref{remarkonnormofwickpowers}.
Here, we used the fact that discrete Wick powers of different orders
are orthogonal. Thus, we even obtain $\lim_{m \rightarrow\infty}\limsup
_{n \rightarrow\infty} \sum_{k = m+1}^{n} \frac{C^{2k}}{k!} = 0$ and,
for all $t \in[0,1]$, a stronger result than condition (\ref{geomwconvstraffganz2}),
\[
\lim_{m \rightarrow\infty} \limsup_{n \rightarrow\infty} \mathbf{E} \Biggl[
\Biggl(\sum_{k = 0}^{n} \frac{a_{n,k}}{k!} (B^{H,n}_{t})^{\diamond k} - \sum
_{k = 0}^{m} \frac{a_{n,k}}{k!}(B^{H,n}_{t})^{\diamond_{n} k} \Biggr)^2 \Biggr] = 0.
\]

By the orthogonality of the different Wick powers, we have
\[
\mathbf{E} \Biggl[ \Biggl(\sum_{k = m+1}^{\infty} \frac
{a_{k}}{k!}(B^{H}_{t})^{\diamond k} \Biggr)^2 \Biggr] = \sum_{k = m+1}^{\infty}
\biggl(\frac{a_{k}}{k!} \biggr)^2 \mathbf{E} [ ((B^{H}_{t})^{\diamond k} )^2 ] \leq
\sum_{k = m+1}^{\infty} \biggl(\frac{C^{2k}}{k!} \biggr) t^{2H k} \rightarrow0
\]
for $m \rightarrow\infty$, which implies that condition (\ref{geomwconvstraffganz3}) even holds in $L^{2}(\Omega,P)$.
\end{pf}

In view of Proposition \ref{propzumdiffconvrate} and Slutsky's
theorem, the previous proposition also implies the following.

\begin{Prop}\label{fdvondiff}
In the context of Theorem \ref{extendedmaintheorem}, convergence
holds in finite-dimensional distributions.
\end{Prop}

\section{Tightness}\label{sectightness}

We now show the tightness of the sequences in Theorems \ref{dattheo1}
and \ref{extendedmaintheorem} by the following criterion, which is a
variant of Theorem 15.6 in Billingsley \cite{Billialt}.

\begin{Theo}\label{LemmaBillsatz}
Suppose that, for the random elements $Y^{n}$ in the Skorokhod space
$D([0,1],\mathbb{R})$ and $\sum_{k=0}^{\infty}\frac{a_{k}}{k!} (B^{H}
)^{\diamond k}$ in $C([0,1],\mathbb{R})$,
\[
Y^{n} \stackrel{fd}{\longrightarrow} \sum_{k=0}^{\infty}\frac
{a_{k}}{k!} (B^{H} )^{\diamond k},
\]
and for $s \leq t$ in $[0,1]$,
\[
\mathbf{E} [ (Y^{n}_{t} - Y^{n}_{s} )^2 ] \leq L \biggl|\frac{ \lfloor nt
\rfloor}{n} - \frac{ \lfloor ns \rfloor}{n} \biggr|^{2H},
\]
where $L > 0$ is a constant. Then $Y^{n}$ converges weakly to $\sum
_{k=0}^{\infty}\frac{a_{k}}{k!} (B^{H} )^{\diamond k}$
in $D([0,1],\mathbb{R})$.
\end{Theo}

\begin{pf}
Let $s<t<u$ in $[0,1]$. By the Cauchy--Schwarz inequality,
\begin{eqnarray*}
&&\mathbf{E} [ |Y^{n}_{t} - Y^{n}_{s} | |Y^{n}_{u} - Y^{n}_{t} | ]\\
&&\quad\leq (\mathbf{E} [ |Y^{n}_{t} - Y^{n}_{s} |^2 ] )^{1/2}
(\mathbf{E} [ |Y^{n}_{u} - Y^{n}_{t} |^2 ] )^{1/2}\\
&&\quad\leq L \biggl|\frac{ \lfloor nt \rfloor}{n} - \frac{ \lfloor ns \rfloor}{n}
\biggr|^{H} \biggl|\frac{ \lfloor nu \rfloor}{n} - \frac{ \lfloor nt \rfloor}{n}
\biggr|^{H} \leq L \biggl|\frac{ \lfloor nu \rfloor}{n} - \frac{ \lfloor ns \rfloor
}{n} \biggr|^{2H}.
\end{eqnarray*}
If $u-s \geq\frac{1}{n}$, we have, since $ \lfloor nu \rfloor\leq nu$
and $- \lfloor ns \rfloor\leq-ns +1$,
\[
\biggl|\frac{ \lfloor nu \rfloor}{n} - \frac{ \lfloor ns \rfloor}{n} \biggr|^{2H}
\leq \bigl(2(u-s) \bigr)^{2H}
\]
and thus
%
\begin{equation}\label{dieglstraffschwexp}
\mathbf{E} [ |Y_{t}^{n} - Y_{s}^{n} | |Y_{u}^{n} - Y_{t}^{n} | ] \leq
L2^{2H} (u-s )^{2H}.
\end{equation}
If $u-s < \frac{1}{n}$, then we have either $ \lfloor ns \rfloor=
\lfloor nt \rfloor$ or $ \lfloor nt \rfloor= \lfloor nu \rfloor$ and
so the left-hand side in (\ref{dieglstraffschwexp}) is zero. Thus,
the inequality (\ref{dieglstraffschwexp}) holds for all $s<t<u$. By
the convergence of the finite-dimensional distributions and \cite{Billialt}, Theorem 15.6, we get the weak convergence of the processes.
\end{pf}

For the application of this criterion to the discrete Wick powers, we
need two lemmas.

\begin{Lemma}\label{dasskalardings}
Let $ (X, \langle\cdot, \cdot\rangle)$ be a real inner product space
and $ \| x \| := \langle x, x \rangle$ the corresponding norm on $X$.
Then, for all $x,y \in X$ and $N \geq1$,
\begin{eqnarray*}
\|x \|^{2N} + \|y \|^{2N} - 2 ( \langle x,y \rangle)^{N}
\leq2^{N+1} ( \|x \| + \|y \| )^{2(N-1)} \|x-y \|^{2}.
\end{eqnarray*}
\end{Lemma}

\begin{pf}
It holds that
\begin{eqnarray*}
2 ( \langle x,y \rangle)^{N} &=&2 \biggl(\frac{1}{2} ( \|x \|^{2} + \|y \|
^{2} - \|x-y \|^{2} ) \biggr)^{N}\\
&=& \frac{1}{2^{N-1}} \Biggl[ ( \|x \|^{2} + \|y \|^{2} )^{N} + \sum_{k =
0}^{N-1} \pmatrix{N\cr k}
 ( \|x \|^{2} + \|y \|^2 )^{k} (-1 )^{N-k} \|x-y \|^{2(N-k)} \Biggr]\\
&=&\frac{1}{2^{N-1}} ( \|x \|^{2} + \|y \|^{2} )^{N}\\
&&{} - \|x-y \|^{2}\frac{1}{2^{N-1}} \sum_{k = 0}^{N-1} \pmatrix{N\cr k}
 ( \|x \|^{2} + \|y \|^2 )^{k} (-1 )^{N-k-1} \|x-y \|^{2(N-k-1)}.
\end{eqnarray*}
Hence, we get
%
\begin{eqnarray}\label{indemtollenlemma2}
&&\|x \|^{2N} + \|y \|^{2N} - 2 ( \langle x,y \rangle)^{N} \nonumber\\
&&\quad= \|x \|
^{2N} + \|y \|^{2N} - \frac{1}{2^{N-1}} ( \|x \|^{2} + \|y \|^{2}
)^{N}\\
&&\qquad{}+ \|x-y \|^{2}\frac{1}{2^{N-1}} \sum_{k = 0}^{N-1} \pmatrix{N\cr k}
 ( \|x \|^{2} + \|y \|^2 )^{k} (-1 )^{N-k-1} \|x-y
 \|^{2(N-k-1)}.\nonumber
\end{eqnarray}

Since $ (\frac{1}{2} )^{N-1} \sum_{k = 0}^{N}{N\choose k} = 2$, the
first line on the right-hand side of (\ref{indemtollenlemma2}) can be
treated as follows:
\begin{eqnarray}\label{nochsonegleichung}
&&\|x \|^{2N} + \|y \|^{2N} - \frac{1}{2^{N-1}} ( \|x \|^{2} + \|y \|
^{2} )^{N}\nonumber\\[-8pt]\\[-8pt]
&&\quad= \frac{1}{2^{N-1}} \sum_{k = 1}^{N-1} \pmatrix{N\cr k}
 \biggl(\frac{ \|x \|^{2N} + \|y \|^{2N}}{2} - \|x \|^{2k} \|y \|^{2(N-k)}
\biggr).\nonumber
\end{eqnarray}
As $ {N\choose k} = {N\choose N-k}$, we now collect the corresponding
summands in sum (\ref{nochsonegleichung}) for $k \neq\frac{N}{2}$. We
obtain, by the mean value theorem with
\[
M_{x,y} := \max_{\lambda\in[0,1]} \bigl(\lambda\|x \| + (1-\lambda) \|y
\| \bigr) = \max\{ \|x \|, \|y \| \}
\]
and since $k(N-k) \leq\frac{N^2}{4}$,
\begin{eqnarray*}
&& \|x \|^{2N} + \|y \|^{2N} - \|x \|^{2k} \|y \|^{2(N-k)} - \|x \|
^{2(N-k)} \|y \|^{2k}\\
&&\quad= ( \|x \|^{2k} - \|y \|^{2k} ) \bigl( \|x \|^{2(N-k)} - \|y \|^{2(N-k)}
\bigr)\\
&&\quad\leq 2k M_{x,y}^{2k-1} \|x-y \| 2(N-k)M_{x,y}^{2(N-k)-1} \|x-y \|
\leq N^2 M_{x,y}^{2(N-1)} \|x-y \|^2.
\end{eqnarray*}
Analogously, we obtain, for $k = \frac{N}{2}$,
\[
\frac{ \|x \|^{2N} + \|y \|^{2N}}{2} - \|x \|^{N} \|y \|^{N} = \frac
{1}{2} ( \|x \|^{N} - \|y \|^{N} )^2 \leq\frac
{1}{2}M_{x,y}^{2(N-1)}N^2 \|x-y \|^2.
\]
Plugging these estimates into (\ref{nochsonegleichung}) and
recalling that $ (\frac{1}{2} )^{N-1} \sum_{k = 1}^{N-1} {N\choose k}
\frac{1}{2} = (1 - \frac{1}{2^{N-1}} )$,
we obtain
%
\begin{equation}\label{dasistdaslemmaundach}
\|x \|^{2N} + \|y \|^{2N} - \biggl(\frac{1}{2} \biggr)^{N-1} ( \|x \|^{2} + \|y \|
^{2} )^{N} \leq \biggl(1 - \frac{1}{2^{N-1}} \biggr)N^2 M_{x,y}^{2(N-1)} \|x-y \|^2.
\end{equation}

For the term in the second line on the right-hand side of (\ref{indemtollenlemma2}), we observe that, by the triangle inequality,
\begin{eqnarray*}
&&( \|x \|^{2} + \|y \|^2 )^{k} (-1 )^{N-k-1} \|x-y \|^{2(N-k-1)}\\
&&\quad\leq ( \|x \| + \|y \| )^{2k} ( \|x \| + \|y \| )^{2(N-k-1)} = ( \|x
\| + \|y \| )^{2(N-1)}.
\end{eqnarray*}
Applying
\[
M_{x,y} \leq \|x \| + \|y \|,\qquad \biggl(\frac{1}{2} \biggr)^{N-1} \sum_{k = 0}^{N-1} \pmatrix{N\cr k} =2-\frac{1}{2^{N-1}},
\]
and (\ref{dasistdaslemmaundach}) to (\ref{indemtollenlemma2}),
we have
\begin{eqnarray*}
&&\|x \|^{2N} + \|y \|^{2N} - 2 ( \langle x,y \rangle)^{N} \\
&&\qquad\leq \biggl[ \biggl(1 - \frac{1}{2^{N-1}} \biggr)N^2 + \biggl(2-\frac{1}{2^{N-1}} \biggr) \biggr] ( \|x \|
+ \|y \| )^{2(N-1)} \|x-y \|^2.
\end{eqnarray*}
By a short calculation and induction, we obtain
\[
\biggl(1 - \frac{1}{2^{N-1}} \biggr)N^2 + \biggl(2-\frac{1}{2^{N-1}} \biggr) \leq\mathbf{1}_{
\{N \neq3 \}} 2^{N} + \mathbf{1}_{ \{N = 3 \}} \frac{17}{2} <2^{N+1}.
\]
\upqed
\end{pf}

\begin{Lemma}\label{Lemmazutight}
For all $t > s$ in $[0,1]$, we have
\[
\frac{1}{N!}\mathbf{E} \bigl[ \bigl((B^{H,n}_{t})^{\diamond_{n} N} -
(B^{H,n}_{s})^{\diamond_{n} N} \bigr)^2 \bigr] \leq 8^N \biggl|\frac{ \lfloor nt
\rfloor}{n} - \frac{ \lfloor ns \rfloor}{n} \biggr|^{2H}.
\]
\end{Lemma}

\begin{pf}
For $N = 1$, the inequality is fulfilled by (\ref{zurboundednessnorm}) in Remark \ref{remarkonboundedness}.
For $N > 1$, we consider the cases $N > \lfloor ns \rfloor$ and $
\lfloor ns \rfloor\geq N$ separately. For $ \lfloor nt \rfloor\geq N
> \lfloor ns \rfloor$, we have $ (B^{H,n}_{s} )^{\diamond_{n} N} = 0$.
Hence, Proposition \ref{walshpropevoll} and Remark \ref{remarkonnormofwickpowers} imply that
\[
\frac{1}{N!}\mathbf{E} \bigl[ \bigl((B^{H,n}_{t})^{\diamond_{n} N} -
(B^{H,n}_{s})^{\diamond_{n} N} \bigr)^2 \bigr] =\frac{1}{N!}\mathbf{E} [
((B^{H,n}_{t})^{\diamond_{n} N} )^2 ] \leq \biggl|\frac{ \lfloor nt \rfloor}{n} \biggr|^{2H N}.
\]
Since $N \geq2$, $2H > 1$ and $\frac{ \lfloor nt \rfloor}{n} \leq1$,
we obtain
\begin{eqnarray*}
\biggl|\frac{ \lfloor nt \rfloor}{n} \biggr|^{2H N} &\leq& \biggl|\frac{ \lfloor nt
\rfloor}{n} \biggr|^{2} = \biggl(\frac{1}{n} \biggr)^2 \bigl( ( \lfloor nt \rfloor- \lfloor
ns \rfloor) + \lfloor ns \rfloor\bigr)^2\\
&\leq& \biggl(\frac{1}{n} \biggr)^2 \bigl( ( \lfloor nt \rfloor- \lfloor ns \rfloor) +
N \bigr)^2 \leq \biggl(\frac{1}{n} \biggr)^2 \bigl( ( \lfloor nt \rfloor- \lfloor ns
\rfloor)(N+1) \bigr)^2\\
&=& (N+1 )^2 \biggl|\frac{ \lfloor nt \rfloor}{n} - \frac{ \lfloor ns \rfloor
}{n} \biggr|^{2}.
\end{eqnarray*}
Since $ (N+1 )^2 \leq3^N$ for $N \geq2$ and $2H < 2$, we obtain
\[
\frac{1}{N!}\mathbf{E} \bigl[ \bigl((B^{H,n}_{t})^{\diamond_{n} N} -
(B^{H,n}_{s})^{\diamond_{n} N} \bigr)^2 \bigr] \leq3^N \biggl|\frac{ \lfloor nt \rfloor
}{n} - \frac{ \lfloor ns \rfloor}{n} \biggr|^{2H}
\]
for all $ \lfloor nt \rfloor\geq N > \lfloor ns \rfloor$. Recall that,
by Proposition \ref{walshpropevoll}, for $ \lfloor nt \rfloor>
\lfloor ns \rfloor\geq N$,
\begin{eqnarray*}
&&\frac{1}{N!}\mathbf{E} \bigl[ \bigl((B^{H,n}_{t})^{\diamond_{n} N} -
(B^{H,n}_{s})^{\diamond_{n} N} \bigr)^2 \bigr] \\
&&\quad\leq\mathbf{E} [(B^{H,n}_{t})^2 ]^{N} + \mathbf{E} [(B^{H,n}_{s})^2
]^{N} - 2\mathbf{E} [(B^{H,n}_{t})(B^{H,n}_{s}) ]^{N}.
\end{eqnarray*}
For any $n \in\mathbb{N}$, we can rewrite
\[
\mathbf{E} [(B^{H,n}_{t})(B^{H,n}_{s}) ] = \sum_{i = 1}^{n}
b_{t,i}^{n}b_{s,i}^{n}
\]
as an ordinary inner product on $\mathbb{R}^{n}$ of the vectors $
(b_{t,1}^{n},\ldots, b_{t,n}^{n} )^{\mathrm{T}}$ and $ (b_{s,1}^{n},\ldots,
b_{s,n}^{n} )^{\mathrm{T}}$. Thus, the application of Lemma \ref{dasskalardings} with $t,s \in[0,1]$ gives
\begin{eqnarray*}
&&\frac{1}{N!}\mathbf{E} \bigl[ \bigl((B^{H,n}_{t})^{\diamond_{n} N} -
(B^{H,n}_{s})^{\diamond_{n} N} \bigr)^2 \bigr]\\
&&\quad\leq 2^{N+1} \bigl(\mathbf{E} [(B^{H,n}_{t})^2 ]^{1/2} + \mathbf
{E} [(B^{H,n}_{s})^2 ]^{1/2} \bigr)^{2(N-1)} \mathbf{E} \bigl[
\bigl((B^{H,n}_{t}) - (B^{H,n}_{s}) \bigr)^2 \bigr]\\
&&\quad\leq 2^{N+1} 2^{2(N-1)} \biggl|\frac{ \lfloor nt \rfloor}{n} - \frac{
\lfloor ns \rfloor}{n} \biggr|^{2H} \leq8^N \biggl|\frac{ \lfloor nt \rfloor}{n} -
\frac{ \lfloor ns \rfloor}{n} \biggr|^{2H}.
\end{eqnarray*}
If $N > \lfloor nt \rfloor$, then the left-hand side of the assertion vanishes.
\end{pf}

\begin{remark}
The proofs for a fractional Brownian motion on some interval $[0,T]
\subset\mathbb{R}$ follow by a straightforward modification. As
$\mathbf{E} [ (B^{H,n}_{t} )^2 ] \leq T^{2H}$ for $t \in[0,T]$ and
\[
\biggl|\frac{ \lfloor nt \rfloor}{n} \biggr|^{2H N} \leq T^{2H(N-1)} \biggl|\frac{
\lfloor nt \rfloor}{n} \biggr|^{2H},
\]
we obtain the previous lemma for $t>s$ in $[0,T]$ as
\[
\frac{1}{N!}\mathbf{E} \bigl[ \bigl((B^{H,n}_{t})^{\diamond_{n} N} -
(B^{H,n}_{s})^{\diamond_{n} N} \bigr)^2 \bigr] \leq (8 T^{2H} )^N \biggl|\frac{
\lfloor nt \rfloor}{n} - \frac{ \lfloor ns \rfloor}{n} \biggr|^{2H}.
\]
\end{remark}

We are now able to prove the weak convergence to the Wick analytic
functionals of a fractional Brownian motion.

\begin{pf*}{Proof of Theorem \ref{dattheo1}}
We apply Theorem \ref{LemmaBillsatz}. The convergence of
finite-dimensional distributions was shown in Proposition \ref{fdmaintheorem}. Let $s<t$ in $[0,1]$. Recall $a_{n,k} \leq C^{k}$. Then, by
the orthogonality of $ ( (B^{H,n}_{t} )^{\diamond_{n} k} - (B^{H,n}_{s}
)^{\diamond_{n} k} )$ for different $k$ and Lemma \ref{Lemmazutight},
we have
\begin{eqnarray*}
&&\mathbf{E} \Biggl[ \Biggl(\sum_{k=0}^{n}\frac{a_{n,k}}{k!} (B^{H,n}_{t}
)^{\diamond_{n} k} - \sum_{k=0}^{n}\frac{a_{n,k}}{k!} (B^{H,n}_{s}
)^{\diamond_{n} k} \Biggr)^2 \Biggr]\\
&&\quad= \sum_{k=0}^{n} \biggl(\frac{a_{n,k}}{k!} \biggr)^2 \mathbf{E} \bigl[
\bigl((B^{H,n}_{t})^{\diamond_{n} k} - (B^{H,n}_{s})^{\diamond_{n} k} \bigr)^2 \bigr]
\leq\sum_{k=0}^{n} \frac{C^{2k}}{k!} 8^{k} \biggl| \frac{ \lfloor nt \rfloor
}{n} - \frac{ \lfloor ns \rfloor}{n} \biggr|^{2H}.
\end{eqnarray*}
Since $0 <\sum _{k=0}^{\infty} \frac{8^{k}C^{2k}}{k!} = \exp(8C^2)=:L
<\infty$, we have
%
\begin{equation}\label{tightwickexp}
\mathbf{E} \Biggl[ \Biggl(\sum_{k=0}^{n}\frac{a_{n,k}}{k!} (B^{H,n}_{t} )^{\diamond
_{n} k} - \sum_{k=0}^{n}\frac{a_{n,k}}{k!} (B^{H,n}_{s} )^{\diamond_{n}
k} \Biggr)^2 \Biggr] \leq L \biggl|\frac{ \lfloor nt \rfloor}{n} - \frac{ \lfloor ns
\rfloor}{n} \biggr|^{2H}.
\end{equation}
\end{pf*}

The alternative approximation, stated in Theorem \ref{extendedmaintheorem}, follows similarly, as we shall now see.

\begin{pf*}{Proof of Theorem \ref{extendedmaintheorem}}
Let $s<t$ in $[0,1]$. Recall that $d_{m,i}^{n} > 0$ only if $i \leq m$.
Thus, by Proposition~\ref{theowalshdiscretewickpowers}, we can write
\[
\sum_{k=0}^{n}\frac{a_{n,k}}{k!} \widetilde{U}^{k,n}_{s} = \sum_{k =
0}^{ \lfloor nt \rfloor} a_{n,k} \mathop{\sum_{C \subseteq\{1,\ldots
, \lfloor nt \rfloor\} }}_{ |C| = k}\mathop{\sum_{m:C \rightarrow\{
1,\ldots, \lfloor ns \rfloor\} }}_{ \mathrm{injective}} \prod_{l \in C}
d_{m(l),l}^n \Xi_{C}^{n}.
\]
Observe that, by the telescoping sum in (\ref{X_teleskopsumme}), we have
\begin{eqnarray*}
\mathop{\mathop{\sum_{m:C \rightarrow\{1,\ldots, \lfloor nt
\rfloor\}}}_{\mathrm{injective}}}_{\exists u : m(u) > \lfloor ns \rfloor} \prod_{l
\in C} d_{m(l),l}^n&\leq& \mathop{\sum_{m:C \rightarrow\{1,\ldots, \lfloor nt \rfloor
\} }}_{ \exists u \dvtx m(u) > \lfloor ns \rfloor} \prod_{l \in C}
d_{m(l),l}^n\\
& =& \sum_{m:C \rightarrow\{1,\ldots, \lfloor nt
\rfloor\} } \prod_{l \in C} d_{m(l),l}^n - \sum_{m:C
\rightarrow\{1,\ldots, \lfloor ns \rfloor\}} \prod_{l \in C}
d_{m(l),l}^n\\
&=& b_{t,C}^{n} - b_{s,C}^{n}.
\end{eqnarray*}
Thus, due to the orthogonality of $\widetilde{U}^{k,n}_{t} - \widetilde
{U}^{k,n}_{s}$ for different values of $k$, Proposition \ref{theowalshdiscretewickpowers} and estimate (\ref{tightwickexp}), we obtain
\begin{eqnarray*}
&&\mathbf{E} \Biggl[ \Biggl(\sum_{k=0}^{n}\frac{a_{n,k}}{k!} \widetilde{U}^{k,n}_{t}
- \sum_{k=0}^{n}\frac{a_{n,k}}{k!} \widetilde{U}^{k,n}_{s} \Biggr)^2 \Biggr]\\
&&\quad= \sum_{k = 0}^{ \lfloor nt \rfloor} a_{n,k}^2 \mathop{\sum_{C
\subseteq\{1,\ldots, \lfloor nt \rfloor\} }}_{ |C| = k} \biggl(\mathop{\mathop{\sum_{
m:C \rightarrow\{1,\ldots, \lfloor nt \rfloor\}
}}_{\mathrm{injective}}}_{\exists u : m(u) > \lfloor ns \rfloor} \prod_{l \in C}
d_{m(l),l}^n \biggr)^2\\
&&\quad\leq\sum_{k = 0}^{ \lfloor nt \rfloor} a_{n,k}^2 \mathop{\sum_{C
\subset\{1,\ldots, \lfloor nt \rfloor\}}}_{ |C| = k} (b_{t,C}^{n} -
b_{s,C}^{n} )^2 \leq L \biggl|\frac{ \lfloor nt \rfloor}{n} - \frac{ \lfloor
ns \rfloor}{n} \biggr|^{2H}
\end{eqnarray*}
and the result follows from Proposition \ref{fdvondiff} and Theorem
\ref{LemmaBillsatz}.
\end{pf*}

\section*{Acknowledgement}

The authors thank Jens-Peter Krei{\ss} for interesting discussions.

\printhistory

\end{document}